\documentclass[11pt,a4paper]{article}
\usepackage[colorlinks]{hyperref}

\usepackage[table]{xcolor}
\usepackage{slashed}
\usepackage{xkeyval}
\usepackage{url,amsmath,ifthen,amssymb,latexsym,pstricks,mathrsfs,comment,amsthm,graphicx,tikz,tikz-cd,enumerate,accents,pgffor,cite,wrapfig,multicol,float}

\usepackage{bibspacing}

\usepackage{enumitem}
\usepackage[margin=10pt,font=small,labelfont=bf, labelsep=period]{caption}

\usepackage[T1]{fontenc}
\usepackage{ wasysym }
\usepackage{stmaryrd}

\usepackage{geometry} 

\begin{document}

\newcommand{\nc}{\newcommand}
\nc{\rnc}{\renewcommand}

\nc\JErev[1]{{\red #1}}

\nc{\noinclude}[1]{#1} 

\nc\ben{\begin{enumerate}[label=\textup{(\roman*)},leftmargin=7mm]}
\nc\bena{\begin{enumerate}[label=\textup{(\alph*)},leftmargin=7mm]}
\nc\een{\end{enumerate}}

\let\oldproofname=\proofname
\rnc{\proofname}{\rm\bf{\oldproofname}}

\nc\lam\lambda
\nc\bn{{\bf n}}
\nc\bu{{\bf u}}
\nc\bv{{\bf v}}
\nc\V{{\bf V}}
\nc\W{{\bf W}}
\nc{\T}{\mathcal T}
\nc{\ze}{\zeta}
\nc{\al}{\alpha}
\nc{\be}{\beta}
\nc{\ga}{\gamma}
\nc{\de}{\delta}
\nc\ka\kappa
\nc\ve\varepsilon
\rnc\th\theta
\nc{\rank}{\operatorname{rank}}
\nc{\Reg}{\operatorname{Reg}}
\nc{\Sub}{\operatorname{Sub}}
\nc{\Eq}{\operatorname{Eq}}
\nc{\im}{\operatorname{im}}
\nc{\Ker}{\operatorname{Ker}}
\nc{\R}{\mathrel{\mathscr R}}
\rnc{\L}{\mathrel{\mathscr L}}
\rnc{\H}{\mathrel{\mathscr H}}
\nc{\D}{\mathrel{\mathscr D}}
\nc{\J}{\mathrel{\mathscr J}}
\rnc{\implies}{\ \Rightarrow\ }
\nc{\sub}{\subseteq}
\nc\bl{{\bf l}}
\nc\br{{\bf r}}
\nc\bk{{\bf k}}
\nc\bq{{\bf q}}
\nc\bs{{\bf s}}
\nc\bt{{\bf t}}
\nc\trans[1]{\left(\begin{smallmatrix}#1\end{smallmatrix}\right)}
\nc\op\oplus
\nc\leqJ{\leq_{\J}}
\nc\geqJ{\geq_{\J}}
\rnc\S{\mathcal S}
\nc\bm{{\bf m}}
\nc\Lam\Lambda
\nc\Om\Omega
\nc\LZ{\operatorname{LZ}}
\nc\RZ{\operatorname{RZ}}
\nc\De\Delta
\nc\RR{\mathbb R}
\nc\E{\mathbb E}
\nc\I{\mathbb I}
\nc\K{\mathbb K}
\nc\bK{{\bf K}}
\nc\Ga\Gamma
\nc\rn\varrho

\nc\bfour{{\bf 4}}
\nc\bthree{{\bf 3}}
\nc\btwo{{\bf 2}}

\nc{\C}{\mathcal C}

\nc{\tfnt}{\lfloor\tfrac n2\rfloor}
\nc{\tfmt}{\lfloor\tfrac m2\rfloor}
\nc{\fnt}{\lfloor\frac n2\rfloor}
\nc{\fnf}{\lfloor\frac n4\rfloor}
\nc{\fmt}{\lfloor\frac m2\rfloor}
\nc{\nt}{\tfrac n2}
\nc{\mmt}{\frac m2}
\nc{\mmmt}{\frac m2}
\nc{\ntf}{\lfloor\frac{n+2}4\rfloor}
\nc{\Divm}{\mathbb D(m)}
\nc{\ind}{\chi_\Z}
\nc{\ds}{\displaystyle}
\nc{\pre}{\preceq}
\nc{\Mod}[1]{\ (\mathrm{mod}\ #1)}
\rnc{\iff}{\ \Leftrightarrow\ }
\nc{\pfcase}[1]{\medskip \noindent {\bf Case #1.}}
\nc{\eitmc}{\end{itemize}\vspace{-2truemm}}
\rnc{\emptyset}{\varnothing}

\rnc{\P}{{\bf P}}
\nc{\Pmn}{\P{m}{n}}
\nc{\Gmn}{G_{m,n}}
\nc{\Bmn}{B_{m,n}}
\nc{\Tmn}{T_{m,n}}

\nc{\fix}{\operatorname{fix}}
\nc{\Fix}{\operatorname{Fix}}
\nc{\la}{\langle}
\nc{\ra}{\rangle}
\nc{\si}{\sigma}
\nc{\sm}{\setminus}
\nc{\set}[2]{\{ {#1} : {#2} \}} 
\nc{\bigset}[2]{\big\{ {#1} : {#2} \big\}} 
\nc{\pf}{\begin{proof}}
\nc{\epf}{\end{proof}}
\nc{\epfres}{\hfill\qed}
\nc{\epfreseq}{\tag*{\qed}}
\nc{\epfeq}{\tag*{\qed}}
\nc{\N}{\mathbb N}
\nc{\Z}{{\bf Z}}
\nc{\mt}{\mapsto}
\nc{\id}{\mathrm{id}}
\nc{\COMMA}{,\qquad}
\nc{\COMMa}{,\ \ \ }
\nc{\OR}{\quad\text{or}\quad}
\nc{\AND}{\qquad\text{and}\qquad}
\nc{\ANd}{\quad\text{and}\quad}
\nc{\ba}{{\bf a}}
\nc{\bb}{{\bf b}}
\nc{\Sum}{\operatorname{\textstyle{\sum}}}

\nc{\bit}{\begin{itemize}}
\nc{\eit}{\end{itemize}}
\nc{\bmc}{\begin{multicols}}
\nc{\emc}{\end{multicols}}
\nc{\itemit}[1]{\item[\emph{(#1)}]}
\nc{\itemnit}[1]{\item[(#1)]}
\nc{\pfitem}[1]{\medskip \noindent #1.}
\nc{\firstpfitem}[1]{#1.}

\numberwithin{equation}{section}

\newtheorem{thm}[equation]{Theorem}
\newtheorem{qu}[equation]{Question}
\newtheorem{lemma}[equation]{Lemma}
\newtheorem{cor}[equation]{Corollary}
\newtheorem{prop}[equation]{Proposition}
\newtheorem{conj}[equation]{Conjecture}
\newtheorem{hope}[equation]{Hope}

\theoremstyle{definition}

\newtheorem{defn}[equation]{Definition}
\newtheorem{rem}[equation]{Remark}
\newtheorem{eg}[equation]{Example}
\newtheorem{prob}[equation]{Problem}

\newtheorem{oprob}{Problem}

\title{\vspace{-3ex}Minimum degrees of finite rectangular bands, null semigroups, and variants of full transformation semigroups\vspace{-5ex}}

\date{}
\author{}

\maketitle

\renewcommand{\thefootnote}{\fnsymbol{footnote}}

\begin{center}
{\large 
Peter J.~Cameron,%
\hspace{-.2em}\footnote{\label{footnote:PJC}Mathematical Institute, University of St Andrews, St Andrews, Fife KY16 9SS, UK. {\it Emails:} {\tt pjc20\,@\,st-andrews.ac.uk}, \ {\tt jdm3\,@\,st-andrews.ac.uk}.}
James East,%
\hspace{-.3em}\footnote{Centre for Research in Mathematics and Data Science, Western Sydney University, Locked Bag 1797, Penrith NSW 2751, Australia. {\it Email:} {\tt j.east\,@\,westernsydney.edu.au}.  Supported by ARC Future Fellowship FT190100632}
Des FitzGerald,%
\hspace{-.3em}\footnote{School of Natural Sciences, University of Tasmania, Private Bag 37, nipaluna/Hobart 7001, Australia.   {\it Email:} {\tt D.FitzGerald\,@\,utas.edu.au}.}
\\
James D.~Mitchell,\hspace{-.2em}\textsuperscript{\ref{footnote:PJC}}
}
Luke Pebody,%
\hspace{-.3em}\footnote{{\it Email:} {\tt luke\,@\,pebody.org}.}
Thomas Quinn-Gregson,%
\hspace{-.3em}\footnote{{\it Email:} {\tt ttquinngregson\,@\,gmail.com}.  This author's research was carried out at TU Dresden with funding from the Deutsche Forschungsgemeinschaft (DFG) and from the European Research Council (Grant Agreement no. 681988, CSP-Infinity).}
\end{center}

\setcounter{footnote}{0}

\renewcommand{\thefootnote}{\arabic{footnote}}

\begin{abstract}
For a positive integer $n$, the full transformation semigroup $\T_n$ consists of all self maps of the set $\{1,\ldots,n\}$ under composition.  Any finite semigroup $S$ embeds in some $\T_n$, and the least such $n$ is called the (minimum transformation) degree of $S$ and denoted $\mu(S)$.  We find degrees for various classes of finite semigroups, including rectangular bands, rectangular groups and null semigroups.  The formulae we give involve natural parameters associated to integer compositions.  Our results on rectangular bands answer a question of Easdown from 1992, and our approach utilises some results of independent interest concerning partitions/colourings of hypergraphs.  

As an application, we prove some results on the degree of a variant $\T_n^a$.  (The variant $S^a=(S,\star)$ of a semigroup $S$, with respect to a fixed element $a\in S$, has underlying set~$S$ and operation $x\star y=xay$.)  It has been previously shown that $n\leq \mu(\T_n^a)\leq 2n-r$ if the sandwich element $a$ has rank~$r$, and the upper bound of $2n-r$ is known to be sharp if $r\geq n-1$.  Here we show that $\mu(\T_n^a)=2n-r$ for $r\geq n-6$.  In stark contrast to this, when $r=1$, and the above inequality says $n\leq\mu(\T_n^a)\leq 2n-1$, we show that $\mu(\T_n^a)/n\to1$ and $\mu(\T_n^a)-n\to\infty$ as $n\to\infty$.  

Among other results, we also classify the $3$-nilpotent subsemigroups of $\T_n$, and calculate the maximum size of such a subsemigroup.

~


~


\end{abstract}

\setcounter{tocdepth}{1}
\tableofcontents


\section{Introduction}\label{sect:intro}

Denote by $\T_X$ the full transformation semigroup over the set $X$, which consists of all self-maps of $X$ under composition.  When $X=\{1,\ldots,n\}$ for some integer $n$ we denote~$\T_X$ by~$\T_n$.  Cayley's Theorem states that every semigroup $S$ embeds in some $\T_X$ with $|X|\leq|S|+1$; see \cite[Theorem~1.1.2]{Howie1995}.  The \emph{(minimum transformation) degree} of a semigroup $S$ is defined to be the cardinal
\[
\mu(S) = {\min}\bigset{|X|}{X\not=\emptyset,\ S \text{ embeds in } \T_X}.
\]
(The requirement that $X\not=\emptyset$ is exclusively to establish the convention that a semigroup of size~$1$ has degree $1$.)
Several authors have calculated $\mu(S)$ for various classes of finite (semi)groups; for semigroups see especially the works of Easdown \cite{Easdown1987,Easdown1988,Easdown1992} and Schein \cite{Schein1992,Schein1988}, but note that in some papers  degrees are defined in terms of representations by \emph{partial} transformations.  (Writing~$\mu'(S)$ for the degree with respect to partial transformations, the reader might like to verify that ${\mu(S)-1\leq\mu'(S)\leq\mu(S)}$.)  See also \cite{BT1991,EEM2017,Holt2010} for computational studies, \cite{Tully1961,Hoehnke1963,Slover1965} for connections to radical theory, and \cite{Wright1975,Johnson1971,EST2010,EH2016,HW2002,BGP1993,KP2000,Saunders2010,Saunders2014} for degrees of finite groups.  The topic is also closely related to enumeration of (classes of) semigroups by size \cite{DM2012,Malandro2019,DK2014,DK2009,Forsythe1955}.
We also mention the very recent study of Margolis and Steinberg \cite{MS2023}, which was written after the current article but published before it, and concerns the class of \emph{Rhodes semisimple semigroups}.  Apart from trivially small exceptions, none of the semigroups we study here belong to this class, so our techniques are necessarily different from theirs.

The initial source of motivation for the current paper came from semigroup variants.  Recall that the \emph{variant} of a semigroup $S$ with respect to a fixed element $a\in S$ is the semigroup $S^a=(S,\star)$, where the \emph{sandwich operation} $\star$ is defined by $x\star y = xay$ for $x,y\in S$.  Variants were introduced by Hickey \cite{Hickey1986,Hickey1983}, building on earlier ideas of Lyapin \cite{Lyapin1960} and Brown \cite{Brown1955}, and have been studied by many others since.  The papers \cite{AM2018,DE2015,Tsyaputa2004,Tsyaputa2005,JE2020,MT2008} study variants of full transformation semigroups; see also \cite{Sandwiches1,Sandwiches2} for a categorical approach.  Sandwich operations also play an important role in computational semigroup theory \cite{EEMP2019}.

Despite the simple definition, the structure of a variant $\T_n^a$ is vastly more complicated than that of $\T_n$ itself.  See for example Figure \ref{fig:egg}, which gives egg-box diagrams for $\T_4$ and a variant~$\T_4^a$, both produced with GAP \cite{GAP,Semigroups}.  (Egg-box diagrams display the structure of a semigroup as determined by Green's relations \cite{Green1951}, which are themselves defined below; see \cite{Howie1995,CPbook} for more details.)  Nevertheless, it was shown in \cite{JE2020} that variants~$\T_n^a$ can be embedded in (ordinary) transformation semigroups of relatively small degree.  For the next statement, which is \cite[Theorem 1.4]{JE2020}, the \emph{rank} of a transformation is the size of its image.

\begin{thm}\label{thm:EM}
If $a\in\T_n$ and $r=\rank(a)$, then $\T_n^a$ embeds in $\T_{2n-r}$.
\end{thm}

In fact, it was shown in \cite{JE2020} that $\T_n^a$ is isomorphic to a \emph{local subsemigroup} of $\T_{2n-r}$ of the form~$b\T_{2n-r}b$, where $b\in\T_{2n-r}$ satisfies $\rank(b)=n$ and $\rank(b^2)=r$.  Theorem \ref{thm:EM} (and the fact that $|\T_n^a|=|\T_n|$) leads to bounds on the degree of a variant:
\begin{equation}\label{eq:2n-r}
n\leq\mu(\T_n^a)\leq 2n-r \qquad\text{for $a\in\T_n$ with $r=\rank(a)$.}
\end{equation}
It was shown in \cite{JE2020} that (for very simple reasons) $\mu(\T_n^a)$ in fact attains the upper bound of $2n-r$ for $r\geq n-1$.  One of our main results, Theorem \ref{thm:n-6} below, improves this by showing that $\mu(\T_n^a)=2n-r$ for $r\geq n-6$.

The situation is very different, however, when the sandwich element $a\in\T_n$ has minimum possible rank $1$, where \eqref{eq:2n-r} becomes $n\leq\mu(\T_n^a)\leq 2n-1$.  For $\rank(a)=1$, we show that:
\bit
\item $\mu(\T_n^a)$ achieves the upper bound of $2n-1$ if and only if $n\leq15$ (Proposition \ref{prop:r=1}),
\item the ratio $\mu(\T_n^a)/n$ tends to $1$ as $n\to\infty$ (Theorem \ref{thm:r=1xn}),
\item the difference $\mu(\T_n^a)-n$ tends to $\infty$ as $n\to\infty$ (Theorem \ref{thm:r=1n+k}).
\eit

En route to proving the above results on variants, we conduct an analysis of various other classes of semigroups.  

Section \ref{sect:RB} concerns rectangular bands.  The main result of this section is Theorem~\ref{thm:RB}, which gives the degree for such a band, extending previously-known results on left and right zero semigroups \cite{GM2008,Easdown1992}.  Table \ref{tab:RB} gives several calculated values.  The proof of Theorem~\ref{thm:RB} utilises new results on hypergraphs, which we believe are of independent interest; see especially Theorem~\ref{thm:H}.  As a further application we also treat rectangular groups in Theorem \ref{thm:RG}.

Section \ref{sect:null} concerns null semigroups, and semigroups we call \emph{right null semigroups} (these include null semigroups and right zero semigroups as special cases).  Among other things, we classify all null subsemigroups of full transformation semigroups and calculate the maximum size of such subsemigroups; see Theorems \ref{thm:null} and \ref{thm:null2}.  Curiously, it transpires that null and left zero semigroups of the same size have the same degree.  We also consider \mbox{($3$-)nilpotent} semigroups in Section \ref{sect:null}, and Theorem \ref{thm:3nil} gives the maximum size of a $3$-nilpotent subsemigroup of $\T_n$; Table \ref{tab:3nil} gives some calculated values.  This result naturally complements the article \cite{DM2012}, which enumerates abstract $3$-nilpotent semigroups by size.  The degree of a uniform right null semigroup is given in Theorem \ref{thm:RN}; calculated values are given in Table \ref{tab:RN}.  We also apply the results of Section~\ref{sect:null} to give an example of a semigroup $S$ such that $\mu(S^a)<\mu(S)$ for all $a\in S$, answering a question from \cite{JE2020}.

The above-mentioned results on variants of $\T_n$ are then given in Section \ref{sect:Tna}.

\begin{figure}[ht]
\begin{center}
\begin{tikzpicture}
\node[above] () at (0,0) {\includegraphics[height=4cm]{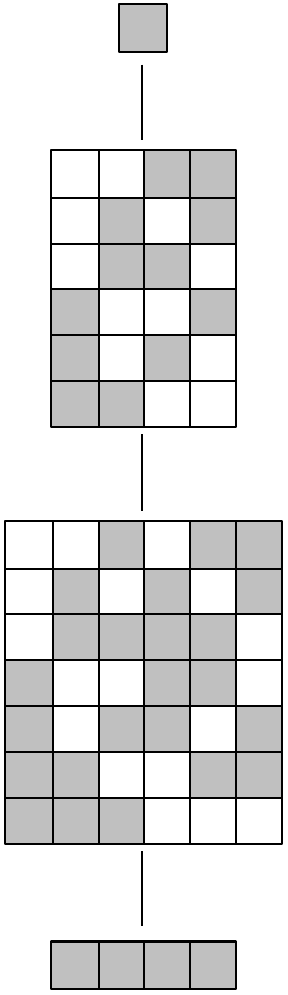}};
\node[above] () at (8,0) {\includegraphics[height=7cm]{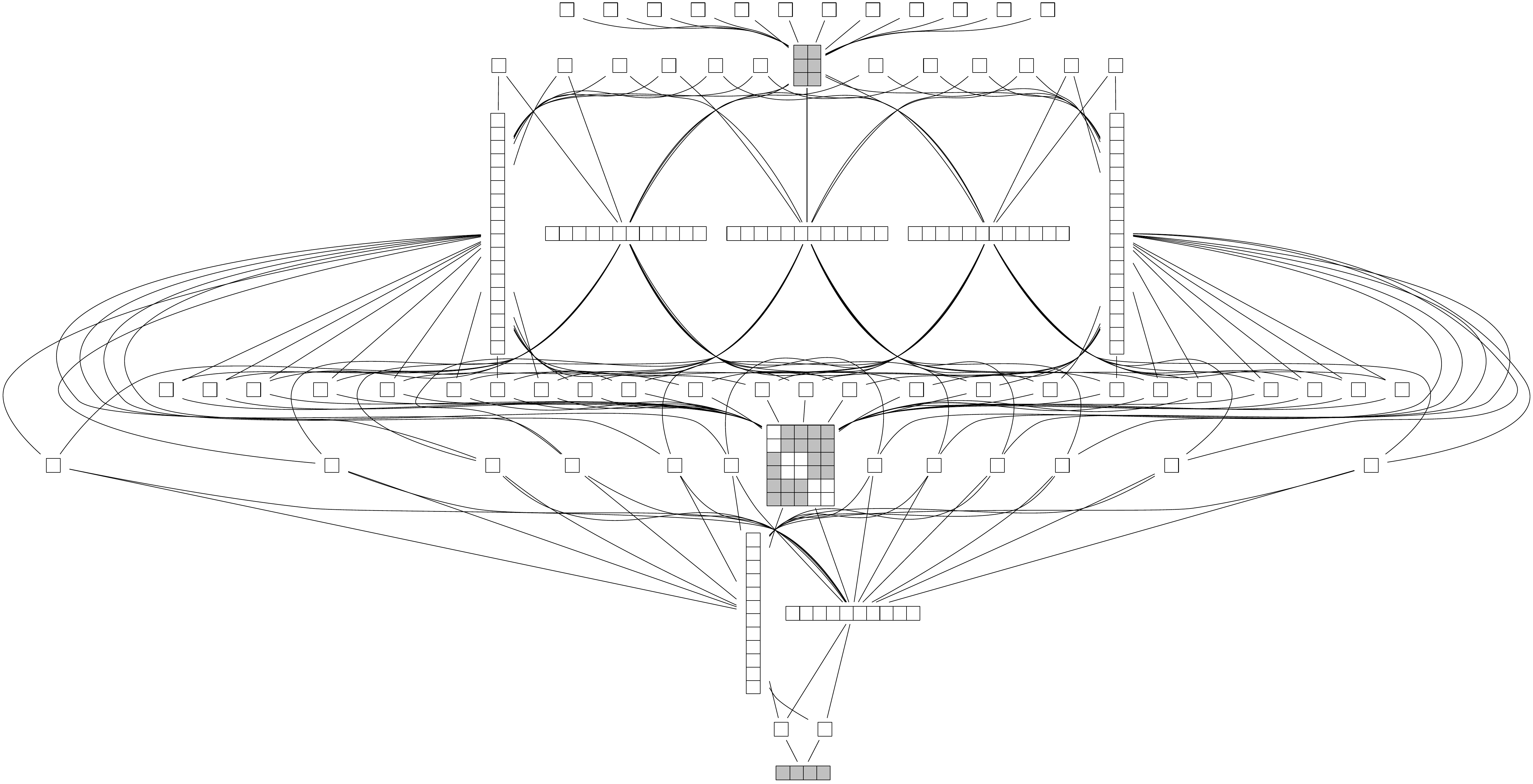}};
\end{tikzpicture}
    \caption[blah]{Egg-box diagrams of the full transformation semigroup $\T_4$ (left) and a variant $\T_4^a$ (right), where $a=\trans{1&2&3&4\\1&2&3&3}$ in standard two-line notation.}
    \label{fig:egg}
   \end{center}
 \end{figure}

\subsection*{Acknowledgements}

We kindly thank the anonymous referees for their careful reading and helpful suggestions.

\section{Preliminaries}\label{sect:prelim}

We begin with some preliminary material on integers and compositions, and (transformation) semigroups.

\subsection{Basic numerical facts}\label{subsect:numbers}

We write $\N=\{1,2,3,\ldots\}$ for the set of natural numbers, and for $n\in\N$ we write $\bn=\{1,\ldots,n\}$.  We adopt the convention that $0^0=1$.

\begin{lemma}\label{lem:fact1}
For integers $n\geq k\geq0$, we have $(n+1)^{k-1}\leq n^k$.
\end{lemma}

\pf
This is clear for $k=0$, and for $k\geq1$ we use the Binomial Theorem:
\[
(n+1)^{k-1} = \sum_{j=0}^{k-1} \tbinom{k-1}jn^{k-1-j}
\leq \sum_{j=0}^{k-1} (k-1)^jn^{k-1-j}
\leq \sum_{j=0}^{k-1} n^jn^{k-1-j}
= k\cdot n^{k-1} \leq n^k.
\qedhere
\]
\epf

\begin{cor}\label{cor:fact1}
For integers $n\geq k\geq0$ and $t\geq0$, we have $(n+t)^{k-t}\leq n^k$.  
\end{cor}

\pf
This is clear for $t=0$ and for $t> k$.  For $1\leq t\leq k$, we repeatedly apply Lemma \ref{lem:fact1}:
\[
n^k \geq (n+1)^{k-1} \geq (n+2)^{k-2} \geq \cdots \geq (n+t)^{k-t}.  \qedhere
\]
\epf

Two closely-related functions $\xi,\al:\N\to\N$ will play an important role throughout, as well as an associated real function $L:\RR^+\to\RR^+$.  Here $\RR^+$ is the set of positive reals.  For $n\in\N$, we define (with the $\bn$ notation introduced at the start of the subsection):
\begin{equation}\label{eq:xial}
\xi(n) = \max\set{t^{n-t}}{t\in\bn} \AND \al(n) = \max\set{t\in\bn}{\xi(n) = t^{n-t}}.
\end{equation}
The numbers $\xi(n)$ and $\al(n)$ appear as Sequences A003320 and A056155 on the OEIS \cite{OEIS}.  Some calculated values are given in Table \ref{tab:xial}.  The $\xi$ and $\al$ functions played a role in \cite{GM2008} in relation to left zero subsemigroups of $\T_n$; they will appear in several other contexts in the current paper.

To define the function $L:\RR^+\to\RR^+$, first consider $f:\RR^+\to\RR$ defined by $f(x)=x+x\ln x$.  Since $f'(x)=2+\ln x$, we see that $f(x)$ is strictly increasing for $x\geq e^{-2}$.  Since $f(e^{-1})=0$, it follows that $f$ restricts to a bijection $\set{x\in\RR}{x>e^{-1}}\to\RR^+$.  We then define $L$ to be the inverse of this bijection.  Thus, for every $t\in\RR^+$, $L(t)$ is the unique real $x>e^{-1}$ satisfying $t=x+x\ln x$.  Note that $L$ is also an increasing function, and that $t\geq L(t) \iff t\geq1 \iff L(t)\geq1$.

\begin{lemma}\label{lem:xi}
\ben
\item \label{xi1}  We have $1=\xi(1)=\xi(2)<\xi(3)<\xi(4)<\cdots$.
\item \label{xi2}  For any $m,n\in\N$, we have $\xi(m)\xi(n)\leq\xi(m+n-1)$.
\item \label{xi3}  For any $n\in\N$, $\al(n)$ is one of $\lfloor x\rfloor$ or $\lceil x\rceil$, where $x=L(n)$, and $\xi(n)=\max\left(\lfloor x\rfloor^{n-\lfloor x\rfloor},\lceil x\rceil^{n-\lceil x\rceil}\right)$.
\een
\end{lemma}

\pf
\firstpfitem{\ref{xi1}}  It is clear that $\xi(1)=\xi(2)=1$.  Now suppose $2\leq m<n$, and write $s=\al(m)$.  Since $m\geq2$, we have $s\geq2$.  Since also $s\in\bm\sub\bn$, we have
\[
\xi(m) = s^{m-s} < s^{n-s} \leq \max\set{t^{n-t}}{t\in\bn} = \xi(n).
\]
\firstpfitem{\ref{xi2}}  Write $s=\al(m)$ and $t=\al(n)$.  Then since $1\leq s+t-1\leq m+n-1$ we have
\begin{align*}
\xi(m+n-1) \geq (s+t-1)^{(m+n-1)-(s+t-1)} &= (s+t-1)^{m-s}(s+t-1)^{n-t} \geq s^{m-s}t^{n-t} = \xi(m)\xi(n).
\end{align*}
\firstpfitem{\ref{xi3}}  As in \cite{GM2008}, define the real function $g:\RR^+\to\RR^+$ by $g(x)=x^{n-x}$.  Then differentiating, we obtain $g'(x)=x^{n-x-1}(n-x-x\ln x)$, so $g'(x)=0$ when $x+x\ln x=n$, i.e.~when $x=L(n)$.  Moreover, if we write $u=L(n)$, then $g(x)$ is increasing for $1\leq x<u$ and decreasing for $x>u$.  The result follows.
\epf

\begin{rem}
It follows that the `max' in the definition of $\al(n)$ in \eqref{eq:xial} is only needed for $n=2$.  Indeed, for $n\geq3$, and writing $u_1=\lfloor u\rfloor$ and $u_2=\lceil u\rceil$, where $u+u\ln u=n$, we have $1<u<n-1$, so that $1\leq u_1<u_2\leq n-1$; it follows that one of $u_1^{n-u_1}$ and $u_2^{n-u_2}$ is even and the other is odd.
\end{rem}

\begin{table}[ht]
\begin{center}
\scalebox{0.9}{
\begin{tabular}{|c|c|r|}
\hline
$n$ & $\al(n)$ & \multicolumn{1}{c|}{$\xi(n)$} \\ 
\hline
\phantom{1}1 & 1 & 1 \\ 
\phantom{1}2 & 2 & 1 \\ 
\phantom{1}3 & 2 & 2 \\ 
\phantom{1}4 & 2 & 4 \\ 
\phantom{1}5 & 3 & 9 \\ 
\phantom{1}6 & 3 & 27 \\ 
\phantom{1}7 & 3 & 81 \\ 
\phantom{1}8 & 4 & 256 \\ 
\phantom{1}9 & 4 & 1024 \\ 
10 & 4 & 4096 \\ 
\hline
11 & 4 & 16384 \\ 
12 & 5 & 78125 \\ 
13 & 5 & 390625 \\ 
14 & 5 & 1953125 \\ 
15 & 6 & 10077696 \\ 
16 & 6 & 60466176 \\ 
17 & 6 & 362797056 \\ 
18 & 6 & 2176782336 \\ 
19 & 7 & 13841287201 \\ 
20 & 7 & 96889010407 \\ 
\hline
\end{tabular}
\qquad\qquad
\begin{tabular}{|c|c|r|}
\hline
$n$ & $\al(n)$ & \multicolumn{1}{c|}{$\xi(n)$} \\ 
\hline
21&  \phantom{1}7&  678223072849\\
22&  \phantom{1}7&  4747561509943\\
23&  \phantom{1}8&  35184372088832\\
24&  \phantom{1}8&  281474976710656\\
25&  \phantom{1}8&  2251799813685248\\
26&  \phantom{1}8&  18014398509481984\\
27&  \phantom{1}9&  150094635296999121\\
28&  \phantom{1}9&  1350851717672992089\\
29&  \phantom{1}9&  12157665459056928801\\
30&  \phantom{1}9&  109418989131512359209\\
\hline
31&  10&  1000000000000000000000\\
32&  10&  10000000000000000000000\\
33&  10&  100000000000000000000000\\
34&  10&  1000000000000000000000000\\
35&  10&  10000000000000000000000000\\
36&  11&  108347059433883722041830251\\
37&  11&  1191817653772720942460132761\\
38&  11&  13109994191499930367061460371\\
39&  11&  144209936106499234037676064081\\
40&  12&  1648446623609512543951043690496\\
\hline
\end{tabular}
}
\caption{Values of the functions $\xi,\al:\N\to\N$ defined in \eqref{eq:xial}.}
\label{tab:xial}
\end{center}
\end{table}

The functions $\xi$ and $\al$ will also come up in another, somewhat indirect, way.  We explore this in the next subsection.

\subsection{Compositions}\label{subsect:comp}

Recall that a \emph{composition} of a natural number $n$ is a tuple $\si=(s_1,\ldots,s_r)$ of positive integers satisfying $s_1+\cdots+s_r=n$, in which case we write $\si\vDash n$ and $|\si|=r$; we call $r$ the \emph{length} of~$\si$. If also $s_1\geq\cdots\geq s_r$, then $\si$ is a \emph{partition} of $n$, and we write $\si\vdash n$.  Certain parameters associated to compositions will play an important role in all that follows, as well as certain numerical functions/sequences defined in terms of them.

Perhaps the simplest parameter associated to a composition $\si=(s_1,\ldots,s_r)$ is its product:
\[
\textstyle{\prod\si} = s_1\cdots s_r.
\]
For integers $1\leq r\leq n$, we define
\begin{equation}\label{eq:pi}
\pi(n) = {\max}\bigset{\textstyle{\prod\si}}{\si\vDash n} \AND \pi_r(n) = {\max}\bigset{\textstyle{\prod\si}}{\si\vDash n,\ |\si|=r}.
\end{equation}
Clearly $\pi(n)=\max_r\pi_r(n)$.  The numbers $\pi(n)$ and $\pi_r(n)$ are well understood, and a proof of the following can be found for example in \cite{GM2008}.  The first part follows from the fact that the maximum value of $\prod\si$, for $\si\vDash n$ with $|\si|=r$, occurs when the entries of $\si$ are `almost equal', in the sense that they are all within $1$ of each other (as $s>t+1$ implies $(s-1)(t+1)= st + s - (t+1) >st$).  The second part then involves showing that the product of almost-equal numbers with a given sum occurs when almost all of the numbers are $3$.

\begin{lemma}\label{lem:pi}
\ben
\item \label{pi1} If $1\leq r\leq n$, then $\pi_r(n) = \lceil \tfrac nr\rceil^t \cdot \lfloor \tfrac nr\rfloor^{r-t}$, where $0\leq t\leq r-1$ is such that $n\equiv t\Mod r$.
\item \label{pi2} For any $n\geq2$ we have
\[
\epfreseq
\pi(n) = \begin{cases}
3^{n/3} &\text{if $n\equiv0\Mod3$}\\
4\cdot3^{(n-4)/3} &\text{if $n\equiv1\Mod3$}\\
2\cdot3^{(n-2)/3} &\text{if $n\equiv2\Mod3$}.
\end{cases}
\]
\een
\end{lemma}

The next family of parameters concerns certain \emph{pairs} of compositions.
Given compositions $\si=(s_1,\ldots,s_r)$ and $\tau=(t_1,\ldots,t_r)$ of the same length (but possibly with different sums), we write $\tau\pre\si$ if $t_i\leq s_i$ for all $i$; for such compositions $\tau\pre\si$, we define
\[
\eta(\si,\tau) = \prod_{i=1}^rt_i^{s_i-t_i}.
\]
Furthermore, for any composition $\si$ we define
\[
\xi(\si) = \max\set{\eta(\si,\tau)}{\tau\pre\si}.
\]
Note that if $\si=(s_1,\ldots,s_r)$, then $\xi(\si) = \xi(s_1)\cdots\xi(s_r)$, which explains our re-use of the $\xi$ symbol.  For $1\leq r\leq n$, we define
\[
\Xi(n) = \max\set{\xi(\si)}{\si\vDash n} \AND \Xi_r(n) = \max\set{\xi(\si)}{\si\vDash n,\ |\si|=r},
\]
so again $\Xi(n)=\max_r\Xi_r(n)$.

\begin{lemma}\label{lem:Xi}
For any $1\leq r\leq n$ we have
\ben
\item \label{Xi1} $\Xi(n) = \xi(n)$,
\item \label{Xi2} $\Xi_r(n) = \xi(n-r+1)$,
\item \label{Xi3} $\Xi_1(n)>\Xi_2(n)>\cdots>\Xi_{n-1}(n)=\Xi_n(n)=1$.
\een
\end{lemma}

\pf
Beginning with \ref{Xi2}, consider a composition $\si=(s_1,\ldots,s_r)\vDash n$.  If $r\geq2$, then from Lemma~\ref{lem:xi}\ref{xi2} and $\xi(1)=1$, we obtain 
\[
\xi(\si) = \xi(s_1)\cdots\xi(s_{r-2})\xi(s_{r-1})\xi(s_r) \leq \xi(s_1)\cdots\xi(s_{r-2})\xi(s_{r-1}+s_r-1)\xi(1) = \xi(\si'),
\]
where $\si'=(s_1,\ldots,s_{r-2},s_{r-1}+s_r-1,1)$.  Continuing, we obtain $\xi(\si)\leq\xi(\tau)$, for the composition $\tau=(n-r+1,1,\ldots,1)$.  Since this is true for all $\si\vDash n$ with $|\si|=r$, it follows that
\[
\Xi_r(n) = \xi(\tau) = \xi(n-r+1)\xi(1)\cdots\xi(1) = \xi(n-r+1).
\]
Now that we have proved \ref{Xi2}, note that \ref{Xi3} then follows from Lemma \ref{lem:xi}\ref{xi1}.  Item \ref{Xi1} quickly follows.
\epf

\subsection{Semigroups}\label{subsect:S}

For more background on semigroups, see \cite{Howie1995,CPbook}.

For a semigroup $S$, we denote by $S^1$ the monoid obtained by adjoining an identity element to $S$ if necessary; so $S=S^1$ if $S$ is a monoid.  Green's $\L$, $\R$ and $\J$ relations are defined for $x,y\in S$ by
\[
x\L y\iff S^1x=S^1y \COMMA
x\R y\iff xS^1=yS^1 \COMMA
x\J y\iff S^1xS^1=S^1yS^1 .
\]
Green's $\H$ and $\D$ relations are defined by ${\H}={\L}\cap{\R}$ and ${\D}={\L}\vee{\R}$, where the latter is the join in the lattice of equivalence relations on $S$.  We have ${\D}={\L}\circ{\R}={\R}\circ{\L}$, and moreover ${\D}={\J}$ if $S$ is finite.

An element $x$ of a semigroup $S$ is \emph{regular} if $x=xyx$ for some $y\in S$ (equivalently, if $x=xyx$ and $y=yxy$ for some $y\in S$), while $x$ is an \emph{idempotent} if $x=x^2$.  We write $\Reg(S)$ and $E(S)$ for the sets of all regular elements and all idempotents of $S$, respectively, neither of which is a subsemigroup in general; both are non-empty when $S$ is finite.  Obviously $E(S)\sub\Reg(S)$. 

An element $x$ of a semigroup $S$ is a \emph{left zero} of $S$ if $xa=x$ for all $a\in S$.  A \emph{left zero semigroup} is a semigroup in which every element is a left zero: i.e., $S$ satisfies the law $xy=x$.  Every left zero semigroup is $\L$-simple, meaning that all elements are $\L$-related.  By a \emph{left zero subsemigroup} of a semigroup $S$ we mean a subsemigroup that happens to be a left zero semigroup (but note that its elements might not be left zeros of $S$ itself).

Right zeros and right zero (sub)semigroups are defined analogously, and they have analogous properties.

If a semigroup $S$ has a left zero and a right zero, then it has a unique left zero and a unique right zero, which are equal, and hence the unique \emph{(two-sided) zero} of $S$.  A \emph{null semigroup} is a semigroup $S$ for which $S^2=\set{xy}{x,y\in S}=\{z\}$ for some $z\in S$, which is then necessarily the zero of $S$.

A \emph{rectangular band} is (isomorphic to) a semigroup of the form $P\times Q$, where $P$ is a left zero semigroup and $Q$ a right zero semigroup; multiplication obeys the rule $(p_1,q_1)(p_2,q_2)=(p_1,q_2)$.  If $|P|=p$ and $|Q|=q$ we say the rectangular band is $p\times q$.  Green's relations on a rectangular band $P\times Q$ are particularly easy to describe:
\[
(p_1,q_1) \L (p_2,q_2) \iff q_1=q_2 \AND
(p_1,q_1) \R (p_2,q_2) \iff p_1=p_2 .
\]
Further, $\H$ is the equality relation, and ${\D}={\J}$ is the universal relation.
As extreme cases, a left zero semigroup of size $p$ is a $p\times1$ rectangular band, with a similar statement for right zero semigroups.

\subsection{Full transformation semigroups}\label{subsect:TX}

The full transformation semigroup $\T_X$ consists of all self-maps of the set $X$, under composition.  When $X$ is the set $\bn=\{1,\ldots,n\}$ for a positive integer $n$, we write $\T_X=\T_n$.  For more on transformation semigroups, see \cite{CPbook,Howie1995,GMbook}.

For $f\in\T_X$ we write $\im(f)=\set{xf}{x\in X}$ for the \emph{image} of $f$, $\rank(f)=|{\im(f)}|$ for the \emph{rank} of $f$, and $\ker(f)=\bigset{(x,y)\in X\times X}{xf=yf}$ for the \emph{kernel} of $f$.  Green's relations on~$\T_X$ are given by
\begin{align*}
f\L g &\iff \im(f)=\im(g) ,\\
f\R g &\iff \ker(f)=\ker(g) ,\\
f\D g\iff f\J f &\iff \rank(f)=\rank(g).
\end{align*}
Fix some subset $\Om$ of $\T_X$.  For any non-empty subset $A\sub X$, any equivalence relation $\ve$ on $X$ (denoted~${\ve\pre X}$), and any cardinal $1\leq r\leq|X|$, we write
\begin{align*}
L_A(\Om) &= \set{f\in\Om}{\im(f)=A},\\
R_\ve(\Om) &= \set{f\in\Om}{\ker(f)=\ve},\\
D_r(\Om) &= \set{f\in\Om}{\rank(f)=r}.
\end{align*}
When there is no chance of confusion, we generally write $L_A=L_A(\T_X)$, and so on.  So the sets $L_A$ ($\emptyset\not=A\sub X$), $R_\ve$ ($\ve\pre X$) and $D_r$ ($1\leq r\leq|X|$) are the $\L$-, $\R$- and ${\J}={\D}$-classes of $\T_X$, respectively.  For an equivalence $\ve\pre X$, we write $\Vert\ve\Vert=|X/\ve|$ for the number of $\ve$-classes.

At times it will be convenient to work with the \emph{kernel partition} of a transformation $f\in\T_X$, defined by $\Ker(f)=X/{\ker(f)}$.  This is a set partition of $X$, and $|{\Ker(f)}| = \Vert{\ker(f)}\Vert = \rank(f)$.

A transformation $f\in\T_X$ will be represented in tabular form as $f=\binom{A_i}{a_i}_{i\in I}$ to indicate that $\im(f)=\set{a_i}{i\in I}$ and $a_if^{-1}=A_i$ for all $i\in I$; we additionally assume that the indexing set $I$ is \emph{faithful} in the sense that $a_i\not=a_j$ whenever $i\not=j$, so that $\rank(f)=|I|$.  Note that $\Ker(f)=\set{A_i}{i\in I}$.  When $\rank(f)$ is finite, we write $f = \trans{A_1&\cdots&A_r\\a_1&\cdots&a_r}$.  It is well known and easy to see that $f=\binom{A_i}{a_i}_{i\in I}$ is idempotent if and only if $a_i\in A_i$ for all $i$.

It will also be convenient to describe here an operation $\op$.  Suppose we have a set of pairwise-disjoint non-empty sets $X_i$ ($i\in I$), and a collection of transformations $f_i\in\T_{X_i}$ ($i\in I$).  Writing $X=\bigcup_{i\in I}X_i$, we define $\bigoplus_{i\in I}f_i$ to be the unique transformation of $X$ whose restriction to each~$X_i$ is $f_i$.  For non-empty subsets $\Om_i\sub\T_{X_i}$ ($i\in I$), we write $\bigoplus_{i\in I}\Om_i$ for the set of all~$\bigoplus_{i\in I}f_i$ with each $f_i\in\Om_i$.  In particular, $\bigoplus_{i\in I}\T_{X_i}$ is a subsemigroup of $\T_X$, and is isomorphic to the (external) direct product $\prod_{i\in I}\T_{X_i}$.

For a fixed non-empty subset $A$ of the set $X$, we have a subsemigroup
\[
\T_X(A) = \set{f\in\T_X}{\im(f)\sub A}.
\]
These subsemigroups have been studied extensively in the literature \cite{SS2008,SS2013,MGS2011,Sanwong2011,Sun2013,FS2014,JE2021}, and will play an important role in the current work.  Note that for any $a\in\T_X$ with $\im(a)=A$, the subsemigroup $\T_X(A)$ is precisely the principal left ideal $\T_Xa=\set{fa}{f\in\T_X}$.  When $X=\bn$ and $A=\bk$ for integers $1\leq k\leq n$, we write $\T_X(A) = \T_n(k)$.

%
While the semigroup $\T_X$ is regular, the same is not true of $\T_X(A)$ in general.  It was shown in \cite{SS2008} that an element $f$ of $\T_X(A)$ is regular (in the semigroup $\T_X(A)$) if and only if $Xf=Af$; an equivalent formulation \cite{JE2021} is that $A$ \emph{saturates} $\ker(f)$, meaning that every $\ker(f)$-class contains an element of $A$.
%

\subsection{Left and right zero transformation semigroups}\label{subsect:LZRZ}

For a subset $\Om$ of a finite semigroup $S$, we write $\lam(\Om)$ and $\rho(\Om)$ for the maximum size of a left or right zero subsemigroup of $S$ contained in $\Om$.  In particular, $\lam(S)$ and $\rho(S)$ are the maximum sizes of a left or right zero subsemigroup of $S$.

Now consider a finite semigroup $S$, and a left zero subsemigroup $T$ of $S$.  Since the elements of~$T$ are all $\L$-related in $T$, it follows that~$T$ is contained in a single $\L$-class of $S$, and hence in a single (regular) $\D$-class; let these classes be $L$ and~$D$, respectively.  Noting that all left zeros are idempotents, we have $T\sub E(L)$.  It is easy to see that $E(L)$ is a left zero semigroup for any regular $\L$-class $L$.  It follows from all this, and the analogous considerations of right zero subsemigroups, that
\begin{align*}
\lam(S) &= \max\set{\lam(D)}{D\in S/{\D}} &&\text{and}& \rho(S) &= \max\set{\rho(D)}{D\in S/{\D}},
\intertext{and that for any $\D$-class $D$ of $S$,}
\lam(D) &= {\max}\bigset{|E(L)|}{L\in D/{\L}} &&\text{and}& \rho(D) &= {\max}\bigset{|E(R)|}{R\in D/{\R}} .
\end{align*}

The paper \cite{GM2008} classified the left and right zero subsemigroups of a finite full transformation semigroup $\T_n$, and calculated the maximum size of such subsemigroups; see also \cite{Easdown1992}.  It will be convenient to interpret these results in our current set-up, using the $\xi$ and $\pi$ parameters from Subsections~\ref{subsect:numbers} and~\ref{subsect:comp}.  For $1\leq r\leq n$, we will write
\[
\lam(n) = \lam(\T_n) \COMMA \rho(n) = \rho(\T_n) \COMMA \lam_r(n) = \lam(D_r(\T_n)) \COMMA \rho_r(n) = \rho(D_r(\T_n)).
\]
As discussed in Subsection \ref{subsect:TX}, the $\L$- and $\R$-classes contained in the $\D$-class $D_r(\T_n)$ are, respectively, the sets of the form
\begin{align*}
L_A = L_A(\T_n) &= \set{f\in\T_n}{\im(f)=A} &&\text{for $A\sub\bn$ with $|A|=r$}\\
R_\ve = R_\ve(\T_n) &= \set{f\in \T_n}{\ker(f)=\ve} &&\text{for $\ve\pre\bn$ with $\Vert\ve\Vert=r$.}
\end{align*}
The size of $E(L_A)$ is $r^{n-r}$, as any idempotent from $L_A$ maps $A$ identically, and maps $\bn\sm A$ arbitrarily into $A$.  The size of $E(R_\ve)$ is equal to the product of the sizes of the $\ve$-classes, as any idempotent from $R_\ve$ maps each $\ve$-class onto a single element of that class.  By maximising these values, one obtains the following, which is \cite[Propositions 3.1 and 3.2]{GM2008}.  

\begin{thm}\label{thm:GM1}
For $1\leq r\leq n$, the maximum size of a left or right zero semigroup contained in $D_r(\T_n)$ is equal to
\[
\epfreseq
\lam_r(n) = r^{n-r} \AND \rho_r(n) = \pi_r(n).
\]
\end{thm}

\begin{rem}\label{rem:GM1}
It is easy to check that for fixed $r$, the sequences $\lam_r(n)$ and $\rho_r(n)$ are non-decreasing in $n\geq r$.  
\end{rem}

\begin{rem}\label{rem:GM1'}
Keeping in mind the formula for $\rho_r(n)=\pi_r(n)$ from Lemma \ref{lem:pi}\ref{pi1}, we of course have $\rho_n(n)=1$ for all $n$.  
Less trivially, we have $\rho_r(n)=2^{n-r}$ for $\frac n2\leq r\leq n$.
\end{rem}

The values of $\lam(n)$ and $\rho(n)$ are found by maximising $\lam_r(n)$ and $\rho_r(n)$ over all $1\leq r\leq n$.  The following is \cite[Theorems 3.3 and 3.4]{GM2008}.  

\begin{thm}\label{thm:GM2}
For $n\geq1$, the maximum size of a left or right zero subsemigroup of  $\T_n$ is equal to
\[
\epfreseq
\lam(n)=\xi(n) \AND \rho(n) = \pi(n).
\]
\end{thm}

We will also need to extend Theorem \ref{thm:GM1} to the semigroups $\T_n(k) = \set{f\in\T_n}{\im(f)\sub\bk}$.  For integers ${1\leq r\leq k\leq n}$, we write 
\[
\lam_r(n,k) = \lam(D_r(\T_n(k))) \AND \rho_r(n,k) = \rho(D_r(\T_n(k)))
\]
for the largest size of a left or right zero semigroup contained in $D_r(\T_n(k))$, respectively.  Note that $D_r(\T_n(k))$ is generally not a $\D$-class of $\T_n(k)$; however, it contains a single regular $\D$-class, as follows from results of \cite[Section 5]{JE2021}.

\begin{lemma}\label{lem:Tnk}
For any $1\leq r\leq k\leq n$, we have $\lam_r(n,k)=\lam_r(n)$ and $\rho_r(n,k)=\rho_r(k)$.
\end{lemma}

\pf
For any $A\sub\bk$ of size $r$, it is clear that $L_A(\T_n)\sub D_r(\T_n(k))$.  Since the $\L$-classes contained in $D_r(\T_n)$ have a common number of idempotents, it quickly follows that ${\lam_r(n,k)=\lam_r(n)}$.

For the statement concerning $\rho_r(n,k)$, consider a regular $\R$-class $R$ contained in $D_r(\T_n(k))$.  Let $\ve$ be the common kernel of the elements of $R$, and let the $\ve$-classes be $A_1,\ldots,A_r$.  Since $R$ is regular, $\bk$ saturates $\ve$, so it follows that $B_i=A_i\cap\bk$ is non-empty for each $i\in\br$; let $\eta\pre\bk$ be the equivalence with classes $B_1,\ldots,B_r$.  Every idempotent from~$R$ is uniquely determined by (and uniquely determines) its restriction to $B_1\cup\cdots\cup B_r=\bk$, which is an idempotent from $R_\eta(\T_k)$; thus, $|E(R)|=|E(R_\eta(\T_k))|$.  Moreover, the above correspondence is reversible; given any $\eta\pre\bk$ with $\Vert\eta\Vert=r$ we can find a suitable $\ve\pre\bn$ by arbitrarily assigning the remaining elements of $\bn\sm\bk$ to $\eta$-classes.  The result now follows.
\epf

One could then obtain results concerning the largest left and right zero subsemigroups of $\T_n(k)$ by maximising the values of $\lam_r(n,k)=\lam_r(n)$ and $\rho_r(n,k)=\rho_r(k)$ over $r\in\bk$.  Note that the latter maximum will be $\rho(k)=\pi(k)$, as in Theorem \ref{thm:GM2}, but the former is ${\max\set{\lam_r(n)}{1\leq r\leq k}}$, and is not necessarily equal to $\lam(n)=\xi(n)$.

\section{Rectangular bands and hypergraphs}\label{sect:RB}

In this section we calculate the (minimum transformation) degree $\mu(B)$ of an arbitrary finite rectangular band $B$, thereby answering a question of Easdown \cite{Easdown1992}.  Since all $p\times q$ rectangular bands are isomorphic, the degree of such a band depends only on the parameters~$p$ and $q$.  Accordingly, for $p,q\in\N$ we define
\[
\be(p,q) = \min\set{n}{\T_n\text{ contains a $p\times q$ rectangular band}}.
\]
Since rectangular bands are $\D$-simple, any such band contained in $\T_n$ is in fact contained in a single $\D$-class $D_r(\T_n)$.  Accordingly, for $p,q,r\in\N$ we define
\[
\be_r(p,q) = \min\set{n}{D_r(\T_n)\text{ contains a $p\times q$ rectangular band}},
\]
and of course $\be(p,q)=\min_r\be_r(p,q)$.  In the above definition of $\be_r(p,q)$, we interpret $\min\emptyset=\infty$.  In particular, we have $\be_1(p,q)=\infty$ for $p\geq2$, since $D_1(\T_n)$ is an $n$-element right zero semigroup (i.e., a $1\times n$ rectangular band) for any $n$.
It follows that $\be(p,q)=\min_{r\geq2}\be_r(p,q)$ when $p\geq2$.

As we noted in Subsection \ref{subsect:S}, the special cases of $q=1$ and $p=1$ correspond to left and right zero semigroups, respectively.  In these cases, values of $\be(p,q)$ and $\be_r(p,q)$ can be quickly deduced from the results of \cite{GM2008,Easdown1992} stated in Theorems \ref{thm:GM1} and \ref{thm:GM2} above.  The following statement uses the~$\xi$ and~$\pi$ parameters from Subsections \ref{subsect:numbers} and \ref{subsect:comp}.

\begin{thm}\label{thm:LZRZ}
\ben
\item \label{LZRZ1} For any $p\geq1$ and any $r\geq2$, we have
\[
\be_r(p,1) = \min\set{n}{r^{n-r}\geq p}=\lceil\log_rp\rceil+r \AND \be(p,1) = \min\set{n}{\xi(n)\geq p}.
\]
\item \label{LZRZ2} For any $q,r\geq1$, we have
\[
\epfreseq
\be_r(1,q) = \min\set{n}{\pi_r(n)\geq q} \AND \be(1,q) = \min\set{n}{\pi(n)\geq q}.
\]
\een
\end{thm}

When $p,q\geq2$, we cannot calculate the numbers $\be(p,q)$ and $\be_r(p,q)$ by finding the largest rectangular bands contained in $\T_n$ and $D_r(\T_n)$.  The reason for this is that even if $\T_n$ contained a rectangular band of size $pq$, it may not contain one with the correct dimentions, $p\times q$.  In fact, it follows from results of \cite[Section 4]{GM2008} that the largest rectangular band contained in $\T_n$ (or in $D_r(\T_n)$) is always a left or right zero semigroup, so such an approach is bound to fail.  Also, since a rectangular band is the direct product of a left and right zero semigroup of appropriate sizes, and since $\mu(S\times T)\leq\mu(S)+\mu(T)$ for any semigroups $S$ and $T$, we have an upper bound of $\be(p,q)\leq\be(p,1)+\be(1,q)$.  This upper bound turns out to be an over-estimate, however, as can be seen by inspecting Table \ref{tab:RB}.

Our approach for $p,q\geq2$ will be via \emph{hypergraphs}.  After proving some results on hypergraphs in Subsection \ref{subsect:H}, which we believe are of independent interest (see especially Theorem~\ref{thm:H}), we return to rectangular bands in Subsection \ref{subsect:RB} (see Theorem \ref{thm:RB}).  We give further applications of our results to rectangular groups in Subsection \ref{subsect:RG} (see Theorem \ref{thm:RG}).

\subsection{Hypergraphs}\label{subsect:H}

A \emph{hypergraph} is a pair $H=(V,\E)$, where $V$ is a set of vertices, and $\E$ is a set of non-empty subsets of $V$ called \emph{(hyper)edges}.  If each edge has size $r$, we say that $H$ is \emph{$r$-uniform} and call it an \emph{$r$-hypergraph}; obviously we must have $r\leq|V|$ if $H$ is non-empty (i.e., has at least one edge).  By an \emph{$r$-partition} of an $r$-hypergraph $H$ we mean a partition $\V=\{V_1,\ldots,V_r\}$ of the vertex set such that each edge is a transversal of $\V$ (i.e., contains a unique point from each block $V_i$ of $\V$).  We call~$H$ \emph{$r$-partite} if such an $r$-partition exists.

Given a partition $\V=\{V_1,\ldots,V_r\}$ of a set $V$, the \emph{complete $r$-partite $r$-hypergraph} $H_\V$ has vertex set $V$, and its edges are all the transversals of $\V$.  The number of edges of $H_\V$ is equal to $|V_1|\cdots|V_r|$, which is the product $\prod\si_\V$ of the composition $\si_\V=(|V_1|,\ldots,|V_r|)$.  Any $r$-hypergraph over $V$ for which $\V$ is an $r$-partition is a subgraph of~$H_\V$, and hence has at most~$\prod\si_\V$ edges.

From the above discussion we have the following well-known result, expressed in terms of the~$\pi$ parameters from Subsection \ref{subsect:comp}.

\begin{lemma}\label{lem:H1}
Let $n,r\geq1$ be integers.  Then there exists a non-empty $r$-partite $r$-hypergraph with~$n$ vertices if and only if $r\leq n$, in which case the maximum number of edges in such a hypergraph is $\pi_r(n)$.  \epfres
\end{lemma}

Our main goal in this subsection is to prove Theorem \ref{thm:H} below, which gives the corresponding result for hypergraphs with some specified number of partitions.  The key technical step in the proof is Proposition \ref{prop:H}, which shows that $r$-hypergraphs with many $r$-partitions give rise to $r$-partite $r$-hypergraphs on smaller vertex sets.  The proof of the proposition is inductive, and relies on the next two lemmas.

\begin{lemma}\label{lem:H2}
Suppose $H$ is an $r$-hypergraph with at least one edge, and at least $rt+1$ distinct $r$-partitions for some integer $t\geq1$.  Then there exists a pair of vertices that belong to the same block in at least $t+1$ of the partitions, but not in all of them.
\end{lemma}

\pf
Let $A=\{v_1,\ldots,v_r\}$ be an edge of $H$, and suppose $\V_1,\ldots,\V_{rt+1}$ are distinct $r$-partitions of $H$.  For each $1\leq i\leq rt+1$, let $\V_i=\{V_{i;1},\ldots,V_{i;r}\}$, where $v_j\in V_{i;j}$ for each $j\in\br$.
Since $\V_1=\{V_{1;1},\ldots,V_{1;r}\}$ and $\V_2=\{V_{2;1},\ldots,V_{2;r}\}$ are distinct, there exists $k\in\br$ such that $V_{1;k}\not\sub V_{2;k}$; fix some $w\in V_{1;k}\sm V_{2;k}$.
For each $i\in\{1,\ldots,rt+1\}$ let $j_i\in\br$ be such that $w\in V_{i;j_i}$.  For each $j\in\br$ let $I_j=\set{i\in\{1,\ldots,rt+1\}}{j_i=j}$.  Since $\{1,\ldots,rt+1\}=I_1\cup\cdots\cup I_r$, we have $|I_m|\geq t+1$ for some $m\in\br$.  
Then $w$ and $v_m$ belong to the same block in at least $t+1$ of the partitions (all the $\V_i$ for $i\in I_m$).

Aiming for a contradiction, suppose $w$ and $v_m$ belong to the same block in \emph{all} of the partitions.  Since $w\in V_{1;k}$ and $v_m \in V_{1;m}$, it follows that $k=m$.  But then $w,v_m\in V_{2;m} = V_{2;k}$, contradicting the definition of $w\in V_{1;k}\sm V_{2;k}$.
\epf

Consider a hypergraph $H$ with vertex set $V$, and let $u,v\in V$ be distinct vertices.   Write $H/\{u,v\}$ for the hypergraph obtained by identifying $u$ and $v$.  Formally, the vertex set of $H/\{u,v\}$ is obtained from $V$ by removing the two vertices $u,v$ and replacing them by a single vertex denoted $uv$; all edges of $H$ involving neither of $u,v$ are still edges of $H/\{u,v\}$; all other edges of $H/\{u,v\}$ are of the form $\{uv,w_1,\dots,w_k\}$ where at least one of $\{u,w_1,\ldots,w_k\}$, $\{v,w_1,\ldots,w_k\}$ or $\{u,v,w_1,\ldots,w_k\}$ is an edge of $H$.  If $H$ is $r$-uniform, then $H/\{u,v\}$ might not be.  Also note that $H/\{u,v\}$ could have fewer edges than $H$, but it cannot have more.

\begin{lemma}\label{lem:H3}
Suppose $H$ is an $r$-hypergraph with $n$ vertices and $q$ edges, with at least~$t+1$ distinct $r$-partitions, for some integer $t\geq1$.  Suppose vertices $u$ and $v$ belong to the same block in $t$ of the $r$-partitions, but not in all of them.  Then $H/\{u,v\}$ is an $r$-hypergraph with $n-1$ vertices and $q$ edges, and has at least $t$ distinct $r$-partitions.
\end{lemma}

\pf
For simplicity, we write $H'=H/\{u,v\}$ throughout the proof.  Let the vertex set of $H$ be $W\cup\{u,v\}$, so the vertex set of $H'$ is $W\cup\{uv\}$.  (This of course has size $n-1$.)

For each edge $A$ of $H$, let $A'$ be the corresponding edge of $H/\{u,v\}$, obtained by replacing~$u$ and/or $v$ by $uv$ if necessary, as explained before the lemma.  To show that $H'$ has $q$ edges, we need to show that the map $A\mt A'$ is injective.  So suppose $A'=B'$ for edges $A,B$ of $H$.  If $A'=B'\sub W$, then $A=A'=B'=B$.  Otherwise, $A'=B'=\{uv\}\cup C$ for some $C\sub W$.  Then~$A$ and $B$ must be one of $\{u\}\cup C$, $\{v\}\cup C$ or $\{u,v\}\cup C$.  But $\{u,v\}\cup C$ cannot be an edge of~$H$, since $u,v$ belong to the same block in some $r$-partition of $H$.  On the other hand, $\{u\}\cup C$ and~$\{v\}\cup C$ cannot both be edges of $H$, since $u,v$ belong to distinct blocks in some $r$-partition of~$H$.  This all shows that $A=B$.

The previous paragraph also shows that $|A'|=|A|=r$ for every edge $A$ of $H$, so $H'$ is $r$-uniform.

Finally, suppose $\V_1,\ldots,\V_t$ are distinct $r$-partitions of $H$ in each of which $u,v$ belong to the same block.  For each $1\leq i\leq t$ let $\V_i=\{V_{i;1},\ldots,V_{i;r}\}$, where $u,v\in V_{i;1}$.  For each such $i$, let $\V_i'=\{V_{i;1}',V_{i;2},\ldots,V_{i;r}\}$, where $V_{i;1}'$ is obtained from $V_{i;1}$ by replacing the pair of vertices $u,v$ (from $H$) by the single vertex $uv$ (from $H'$).  Then by the above characterisation of the edges of~$H$, it is clear that each $\V_i'$ is an $r$-partition of $H'$, and that these are still distinct.
\epf

\begin{prop}\label{prop:H}
If $H$ is an $r$-hypergraph with $n$ vertices and $q\geq1$ edges, with at least $r^{l-1}+1$ distinct $r$-partitions for some integer $l\geq1$, then there is an $r$-partite $r$-hypergraph with $n-l$ vertices and $q$ edges.
\end{prop}

\pf
The proof is by induction on $l$.  Suppose first that $l=1$.  So $H$ has at least two distinct $r$-partitions, say $\V_1$ and $\V_2$.  Since these are distinct, there exist vertices $u,v$ belonging to the same block in $\V_1$ but not in $\V_2$.  By Lemma \ref{lem:H3} (with $t=1$), $H/\{u,v\}$ is an $r$-partite $r$-hypergraph with $n-1=n-l$ vertices and $q$ edges.

Now suppose $l\geq2$.  Since $r^{l-1} +1 = r\cdot r^{l-2} + 1$, it follows from Lemma \ref{lem:H2} (with $t=r^{l-2}$) that there are vertices $u,v$ of $H$ that belong to the same block in at least $r^{l-2}+1$ of the partitions, but not in all of them.  By Lemma \ref{lem:H3} (with $t=r^{l-2}+1$), $H/\{u,v\}$ is an $r$-hypergraph with $n-1$ vertices and $q$ edges, and with at least $r^{l-2}+1=r^{(l-1)-1}+1$ distinct $r$-partitions.  By induction there exists an $r$-partite $r$-hypergraph with $(n-1)-(l-1)=n-l$ vertices and $q$ edges.
\epf

Here is the main result of this subsection.

\begin{thm}\label{thm:H}
Let $n,p\geq1$ and $r\geq2$ be integers, and let $l=\lceil\log_rp\rceil$.  Then there exists a non-empty $r$-hypergraph with $n$ vertices and at least $p$ distinct $r$-partitions if and only if $r+l\leq n$, in which case the maximum number of edges in such a hypergraph is $\pi_r(n-l)$.
\end{thm}

\pf
If $p=1$ then $l=0$, and the result reduces to Lemma \ref{lem:H1}.  For the rest of the proof we assume that $p\geq2$.  It follows from the definition of $l$ that $l\geq1$, and $r^{l-1}+1\leq p\leq r^l$.

Suppose there exists a non-empty $r$-hypergraph with $n$ vertices and $q(\geq1)$ edges, and with~$p$ distinct $r$-partitions.  Since $p\geq r^{l-1}+1$, Proposition \ref{prop:H} applies, and it tells us there is an $r$-partite $r$-hypergraph with $n-l$ vertices and $q$ edges.  It then follows from Lemma~\ref{lem:H1} that $r\leq n-l$ (i.e., $r+l\leq n$) and that $q\leq\pi_r(n-l)$.

It remains to show that if $r+l\leq n$ there exists an $r$-hypergraph with $n$ vertices and $\pi_r(n-l)$ edges, with (at least) $p$ distinct $r$-partitions.  To do so, suppose $r+l\leq n$.  By Lemma~\ref{lem:H1} we may fix an $r$-hypergraph $H$ on $n-l$ vertices with $\pi_r(n-l)$ edges.  Let $\V=\{V_1,\ldots,V_r\}$ be an $r$-partition of $H$.  Now let $H'$ be the $r$-hypergraph obtained from $H$ by adding $l$ isolated vertices.  Then $H'$ has $n$ vertices and $\pi_r(n-l)$ edges, and (at least) $r^l$ distinct $r$-partitions (each of the new vertices, independently, can go in any of the blocks $V_1,\ldots,V_r$).  Since $r^l\geq p$, we are done.
\epf

\subsection{Rectangular bands}\label{subsect:RB}

We now use the results of the previous subsection to calculate the degree of an arbitrary finite rectangular band, answering a question of Easdown \cite{Easdown1992}.  The first result shows how to construct certain `large' such bands.  Its proof shares some ideas with that of Theorem \ref{thm:H}.

\begin{prop}\label{prop:RB}
If $1\leq r\leq n$, and if $0\leq l\leq n-r$, then $D_r(\T_n)$ contains an $r^l\times\pi_r(n-l)$ rectangular band.
\end{prop}

\pf
Fix a partition $\V=\{V_1,\ldots,V_r\}$ of $\bn\sm\bl=\{l+1,\ldots,n\}$ with $|V_1|\cdots|V_r|=\pi_r(n-l)$.  For each function $f:\bl\to\br$, and each transversal $A=\{v_1,\ldots,v_r\}$ of $\V$, with each $v_i\in V_i$, we define a transformation
\[
e(f,A) = \trans{V_1\cup1f^{-1}&\cdots&V_r\cup rf^{-1}\\v_1&\cdots&v_r} \in D_r(\T_n).
\]
It is easy to check that these compose according to the rule $e(f,A)\cdot e(g,B) = e(f,B)$.  It follows that the set of all such $e(f,A)$ is a rectangular band contained in $D_r(\T_n)$.  Since there are $r^l$ functions $\bl\to\br$, and $\pi_r(n-l)$ transversals of $\V$, the band has the stated dimensions.
\epf

\begin{lemma}\label{lem:RB2}
If $D_r(\T_n)$ contains a $p\times q$ rectangular band, where $p,q\geq2$, then with $l=\lceil\log_rp\rceil$ we have $\pi_r(n-l)\geq q$.
\end{lemma}

\pf
Suppose $B$ is a $p\times q$ rectangular band contained in $D_r(\T_n)$.  Since $D_1(\T_n)$ is a right zero semigroup, and since $p\geq2$, we have $r\geq2$.  Let
\[
\I = \set{\im(f)}{f\in B} \AND \K = \set{\Ker(f)}{f\in B}
\]
be the sets of all images and kernel partitions of the elements of $B$.  

Since each element of $B$ has rank $r$, each member of $\I$ is a subset of $\bn$ of size $r$, which means that $H=(\bn,\I)$ is an $r$-hypergraph.  
Moreover, given $A\in\I$ and $\bK\in\K$ we have $A=\im(f)$ and $\bK=\Ker(g)$ for some $f,g\in B$, and then since $fg\in B\sub D_r$ it follows that~$A$ is a transversal of~$\bK$.  This means that each member of $\K$ is an $r$-partition of $H$.

Since $B$ has $p$ $\R$-classes and $q$ $\L$-classes, it follows that $|\I|=q$ and $|\K|=p$.  This all means that $H$ is an $r$-hypergraph with $n$ vertices and $q$ edges, and with $p$ distinct $r$-partitions.  It then follows from Theorem \ref{thm:H} that $\pi_r(n-l)\geq q$.
\epf

We can now give an expression for $\be(p,q)$, which we recall denotes the degree of a $p\times q$ rectangular band.  The cases where one or both of $p,q$ is $1$ were treated in Theorem \ref{thm:LZRZ}.

\begin{thm}\label{thm:RB}
If $p,q\geq2$, then the degree of a $p\times q$ rectangular band is equal to
\[
\be(p,q)=\min_{r\geq2}\be_r(p,q), \qquad\text{where}\qquad \be_r(p,q) = {\min}\bigset{n}{\pi_r(n-\lceil\log_rp\rceil)\geq q} \text{ for $r\geq2$.}
\]
\end{thm}

\pf
At the beginning of Section \ref{sect:RB} we explained that $\be(p,q)=\min_{r\geq2}\be_r(p,q)$, where
\[
\be_r(p,q) = \min\set{n}{D_r(\T_n)\text{ contains a $p\times q$ rectangular band}}.
\]
It therefore remains to show that for $n,r\geq2$, $D_r(\T_n)$ contains a $p\times q$ rectangular band if and only if $\pi_r(n-l)\geq q$, where $l=\lceil\log_rp\rceil$.  The forwards implication is Lemma \ref{lem:RB2}.  

For the converse, suppose $\pi_r(n-l)\geq q$.  Then of course we have $n-l\geq r$, and so $l\leq n-r$.  Proposition \ref{prop:RB} then tells us that $D_r(\T_n)$ contains an $r^l\times\pi_r(n-l)$ rectangular band, and hence also a $p\times q$ rectangular band, since $p\leq r^l$ and $q\leq\pi_r(n-l)$.
\epf

\begin{rem}\label{rem:RB}
In principle, the expression for $\be(p,q)$ in Theorem \ref{thm:RB} involves finding the minimum of an infinite sequence of numbers: i.e., the $\be_r(p,q)$ for all $r\geq2$.  In practice however, only a finite subsequence needs to be considered.  For example, to find $\be(p,q)$, one could first find (say) $m=\be_2(p,q)$.  We then know that $\be(p,q)\leq m$, and since $\be_r(p,q)\geq r$ for all $r$ it quickly follows that $\be(p,q)=\min_{2\leq r\leq m}\be_r(p,q)$.

Table \ref{tab:RB} gives several values of $\be(p,q)$, calculated in this way.  The table also gives the ranks~$r$ such that $\be(p,q)=\be_r(p,q)$.  For some values of $p,q$ there are multiple such $r$, in which case the table gives the minimum $r$.
\end{rem}

\begin{rem}
A \emph{band} is a semigroup consisting entirely of idempotents.  It is well known that every band is a semilattice of rectangular bands; see \cite[Section 4.4]{Howie1995} for a proof of this fact, and an explanation of the terminology.  We leave it as an open problem to determine whether this structure (and Theorem \ref{thm:RB}) leads to a formula for the degree of an arbitrary finite band in terms of its underlying semilattice and constituant rectangular bands.
The degree of a semilattice of groups was determined by Easdown in \cite{Easdown1988}.
\end{rem}

\begin{table}[htp]
\begin{center}
\begin{tabular}{|l|cccccccccc|}
\hline
$	p\sm q	$ & $	1	$ & $	2	$ & $	3	$ & $	4	$ & $	5	$ & $	6	$ & $	7	$ & $	8	$ & $	9	$ & $	10	$ \\
\hline
$\phantom{1}1$ & 1, 1 & 2, 1 & 3, 1 & 4, 1 & 5, 1 & 5, 2 & 6, 2 & 6, 2 & \phantom{1}6, 2 & \phantom{1}7, 2 \\
$\phantom{1}2$ & 3, 2 & 4, 2 & 5, 2 & 5, 2 & 6, 2 & 6, 2 & 7, 2 & 7, 2 & \phantom{1}7, 2 & \phantom{1}8, 2 \\
$\phantom{1}3$ & 4, 2 & 5, 2 & 6, 2 & 6, 2 & 7, 2 & 7, 2 & 7, 3 & 7, 3 & \phantom{1}8, 2 & \phantom{1}8, 3 \\
$\phantom{1}4$ & 4, 2 & 5, 2 & 6, 2 & 6, 2 & 7, 2 & 7, 2 & 8, 2 & 8, 2 & \phantom{1}8, 2 & \phantom{1}9, 2 \\
$\phantom{1}5$ & 5, 2 & 6, 2 & 7, 2 & 7, 2 & 8, 2 & 8, 2 & 8, 3 & 8, 3 & \phantom{1}9, 2 & \phantom{1}9, 3 \\
$\phantom{1}6$ & 5, 2 & 6, 2 & 7, 2 & 7, 2 & 8, 2 & 8, 2 & 8, 3 & 8, 3 & \phantom{1}9, 2 & \phantom{1}9, 3 \\
$\phantom{1}7$ & 5, 2 & 6, 2 & 7, 2 & 7, 2 & 8, 2 & 8, 2 & 8, 3 & 8, 3 & \phantom{1}9, 2 & \phantom{1}9, 3 \\
$\phantom{1}8$ & 5, 2 & 6, 2 & 7, 2 & 7, 2 & 8, 2 & 8, 2 & 8, 3 & 8, 3 & \phantom{1}9, 2 & \phantom{1}9, 3 \\
$\phantom{1}9$ & 5, 3 & 6, 3 & 7, 3 & 7, 3 & 8, 3 & 8, 3 & 8, 3 & 8, 3 & \phantom{1}9, 3 & \phantom{1}9, 3 \\
$10$ & 6, 2 & 7, 2 & 8, 2 & 8, 2 & 9, 2 & 9, 2 & 9, 3 & 9, 3 & 10, 2 & 10, 3 \\
\hline
\end{tabular}
\\[5mm]
\begin{tabular}{|l|cccccccccc|}
\hline
$	\phantom{1}p\sm q	$ & $	10	$ & $	20	$ & $	30	$ & $	40	$ & $	50	$ & $	60	$ & $	70	$ & $	80	$ & $	90	$ & $	100	$ \\
\hline
$\phantom{1}10$ & 10, 3 & 11, 4 & 12, 4 & 13, 4 & 13, 4 & 14, 4 & 14, 4 & 14, 4 & 15, 4 & 15, 4 \\
$\phantom{1}20$ & 10, 3 & 12, 3 & 12, 5 & 13, 5 & 14, 4 & 14, 5 & 14, 5 & 15, 4 & 15, 5 & 15, 5 \\
$\phantom{1}30$ & 11, 3 & 12, 4 & 13, 4 & 14, 4 & 14, 4 & 14, 6 & 15, 4 & 15, 4 & 15, 6 & 16, 4 \\
$\phantom{1}40$ & 11, 3 & 12, 4 & 13, 4 & 14, 4 & 14, 4 & 15, 4 & 15, 4 & 15, 4 & 16, 4 & 16, 4 \\
$\phantom{1}50$ & 11, 3 & 12, 4 & 13, 4 & 14, 4 & 14, 4 & 15, 4 & 15, 4 & 15, 4 & 16, 4 & 16, 4 \\
$\phantom{1}60$ & 11, 3 & 12, 4 & 13, 4 & 14, 4 & 14, 4 & 15, 4 & 15, 4 & 15, 4 & 16, 4 & 16, 4 \\
$\phantom{1}70$ & 11, 3 & 13, 3 & 13, 5 & 14, 5 & 15, 4 & 15, 5 & 15, 5 & 16, 4 & 16, 5 & 16, 5 \\
$\phantom{1}80$ & 11, 3 & 13, 3 & 13, 5 & 14, 5 & 15, 4 & 15, 5 & 15, 5 & 16, 4 & 16, 5 & 16, 5 \\
$\phantom{1}90$ & 12, 3 & 13, 4 & 13, 5 & 14, 5 & 15, 4 & 15, 5 & 15, 5 & 16, 4 & 16, 5 & 16, 5 \\
$100$ & 12, 3 & 13, 4 & 13, 5 & 14, 5 & 15, 4 & 15, 5 & 15, 5 & 16, 4 & 16, 5 & 16, 5 \\
\hline
\end{tabular}
\\[5mm]
\begin{tabular}{|l|cccccccccc|}
\hline
$	p\sm q	$ & $	10^0	$ & $	10^1	$ & $	10^2	$ & $	10^3	$ & $	10^4	$ & $	10^5	$ & $10^6	$ & $	10^7	$ & $	10^8	$ & $10^9$	 \\
\hline
$	10^0	$ &  \phantom{1}1, 1 & \phantom{1}7, 2 & 13, 4 & 20, 5 & 26, 8 & 32, 10 & 38, 13 & 45, 14 & 51, 16 & 57, 19 \\
$	10^1	$ & \phantom{1}6, 2 & 10, 3 & 15, 4 & 21, 10 & 27, 10 & 33, 10 & 39, 13 & 46, 14 & 52, 16 & 58, 19 \\
$	10^2	$ & \phantom{1}8, 3 & 12, 3 & 16, 5 & 22, 10 & 28, 10 & 34, 10 & 40, 13 & 47, 14 & 53, 16 & 59, 19 \\
$	10^3	$ & \phantom{1}9, 4 & 13, 4 & 18, 4 & 23, 10 & 29, 10 & 35, 10 & 41, 13 & 48, 14 & 54, 16 & 60, 19 \\
$	10^4	$ & 11, 4 & 15, 4 & 19, 5 & 24, 10 & 30, 10 & 36, 10 & 42, 13 & 49, 14 & 55, 16 & 61, 19 \\
$	10^5	$ & 13, 4 & 17, 4 & 20, 7 & 25, 10 & 31, 10 & 37, 10 & 43, 13 & 49, 18 & 55, 18 & 61, 19 \\
$	10^6	$ & 14, 4 & 18, 4 & 22, 5 & 26, 10 & 32, 10 & 38, 10 & 44, 13 & 50, 16 & 56, 16 & 62, 19 \\
$	10^7	$ & 15, 6 & 19, 6 & 23, 6 & 27, 10 & 33, 10 & 39, 10 & 45, 13 & 51, 15 & 57, 16 & 63, 19 \\
$	10^8	$ & 17, 5 & 21, 5 & 24, 7 & 28, 10 & 34, 10 & 40, 10 & 46, 13 & 52, 14 & 58, 16 & 64, 19 \\
$	10^9	$ & 18, 5 & 22, 5 & 25, 7 & 29, 10 & 35, 10 & 41, 10 & 47, 13 & 53, 14 & 59, 16 & 64, 20 \\
\hline
\end{tabular}
\caption{Calculated values of $\be(p,q)$.  If the $(p,q)$ entry is $n,r$, then $n=\be(p,q)=\be_r(p,q)$, with $r$ minimal.  This means that $n$ is minimal with the property that $\T_n$ contains a $p\times q$ rectangular band, and $1\leq r\leq n$ is minimal with the property that $D_r(\T_n)$ contains such a band.}
\label{tab:RB}
\end{center}
\end{table}

\subsection{Rectangular groups}\label{subsect:RG}

Recall that a \emph{rectangular group} is (isomorphic to) a direct product $S=B\times G$, where $B$ is a rectangular band and $G$ is a group.  If the band $B$ is $p\times q$, we say that $S$ is a \emph{$p\times q$ rectangular group over $G$}.  An alternative characterisation of rectangular groups can be found in \cite[Exercise~10, page~139]{Howie1995}: they are precisely the $\D$-simple regular semigroups whose idempotents form a subsemigroup.

The results of previous subsections can be applied to calculate the degree of an arbitrary finite rectangular group:

\begin{thm}\label{thm:RG}
If $S$ is a finite $p\times q$ rectangular group over $G$, then 
\[
\mu(S) = \min_{r\geq\mu(G)}\be_r(p,q).
\]
The numbers $\be_r(p,q)$ are given in Theorems \ref{thm:LZRZ} and \ref{thm:RB}.
\end{thm}

\pf
Consider a finite $p\times q$ rectangular group $S=B\times G$, and for simplicity write
\[
n=\mu(S) \AND N=\min_{r\geq\mu(G)}\be_r(p,q).
\]
We must show that $n=N$.

By definition of $n$, there exists an embedding $\phi:S\to\T_n$.  Since $S$ is $\D$-simple, $\im(S)$ is contained in a single $\D$-class of $\T_n$, say $D_r(\T_n)$.  Since $G$ embeds in~$S$, it also embeds into a group $\H$-class in $D_r(\T_n)$; since such an $\H$-class is isomorphic to the symmetric group $\S_r$, it follows that $r\geq\mu(G)$.  Similarly, $D_r(\T_n)$ contains an isomorphic copy of $B$, and so $n\geq\be_r(p,q)$ by definition.  This all shows that $n\geq N$.

To show that $n\leq N$, we must show that $S$ embeds in $\T_N$.  To do so, suppose $r\geq\mu(G)$ is such that $N=\be_r(p,q)$.  By definition, $D_r(\T_N)$ contains a $p\times q$ rectangular band $C$.  Let $T=\bigcup_{c\in C}H_c$ be the union of all the $\H$-classes of $\T_N$ containing idempotents from $C$, so that $T$ is a $p\times q$ rectangular group over $\S_r$ (this follows quickly from the above-mentioned characterisation of rectangular groups from \cite{Howie1995}).  Since $G$ embeds in $\S_r$ (as $r\geq\mu(G)$), it follows that $S=B\times G$ embeds in $B\times\S_r\cong T$, and hence in $\T_N$.
\epf

\begin{rem}
As in Remark \ref{rem:RB}, the expression for $\mu(S)$ in Theorem \ref{thm:RG} involves finding the minimum of an infinite sequence of numbers.  However, we again have $\mu(S) = \min_{\mu(G)\leq r\leq m}\be_r(p,q)$ for any \emph{a priori} known upper bound $m$ of $\mu(S)$.  For example, we can take $m=\be(p,q)+\mu(G)$, since $S$ is the direct product of a $p\times q$ rectangular band with $G$.
\end{rem}

\section{(Right) null and nilpotent semigroups}\label{sect:null}

In this section we study null subsemigroups (Subsections \ref{subsect:null} and \ref{subsect:maxnull}) of a full transformation semigroup~$\T_X$, as well as some related kinds of subsemigroups: nilpotent subsemigroups (Subsections~\ref{subsect:nil} and \ref{subsect:maxnil}), and so-called right null subsemigroups (Subsection~\ref{subsect:RN}).  As an application of our results on null semigroups, we give an example of a semigroup with a larger degree than one of its variants (Subsection \ref{subsect:Sa<S}).  Our results on right null semigroups will be used in Section~\ref{sect:Tna}.

\subsection{Null semigroups}\label{subsect:null}

Let $X$ be an arbitrary set, possibly infinite.  Obviously the zero of a null subsemigroup of $\T_X$ is an idempotent.  So we fix an idempotent ${\ze\in E(\T_X)}$, with the goal of classifying the null subsemigroups containing~$\ze$.  We write ${\ze=\binom{Z_i}{z_i}_{i\in I}}$, noting that $z_i\in Z_i$ for all $i$, as $\ze$ is an idempotent.  

\begin{defn}\label{defn:ze}
Given an idempotent $\ze=\binom{Z_i}{z_i}_{i\in I}$ of $\T_X$, a \emph{$\ze$-system} is a collection $\W$ of sets
\[
\W = \set{W_i}{i\in I} \qquad\text{such that}\qquad z_i\in W_i \sub Z_i \text{ for all $i\in I$.}
\]
For such a $\ze$-system $\W$, we define
\[
N(\W,\ze) = \set{f\in\T_X}{Z_if\sub W_i \text{ and } W_if=\{z_i\} \text{ for all $i\in I$}}.
\]
\end{defn}

\begin{lemma}\label{lem:NWz}
Given an idempotent $\ze$ of $\T_X$, a subset $S$ of $\T_X$ is a null semigroup with zero $\ze$ if and only if $\ze\in S\sub N(\W,\ze)$ for some $\ze$-system $\W$.
\end{lemma}

\pf
\firstpfitem{($\Rightarrow$)}  Suppose $S$ is a null semigroup with zero $\ze$.  Obviously $\ze\in S$.  If $f\in S$, then from $\ze=f\ze$ we obtain $Z_i f\sub Z_i$ for all $i$.  Since also $\{z_i\}=Z_i\ze$ for all $i$, it quickly follows that $\W=\set{W_i}{i\in I}$ is a $\ze$-system, where
\[
W_i = \bigcup_{f\in S}Z_if \qquad\text{for all $i\in I$}.
\]
To see that $S\sub N(\W,\ze)$, let $f\in S$ and $i\in I$ be arbitrary.  We have $Z_if\sub W_i$ by definition.  Moreover, for any $g\in S$ we have $\{z_i\}=Z_i\ze=(Z_ig)f$, from which it follows that $W_if=\{z_i\}$.

\pfitem{($\Leftarrow$)}  Now suppose $\ze\in S\sub N(\W,\ze)$ for some $\ze$-system $\W=\set{W_i}{i\in I}$.  Let $f,g\in S$ be arbitrary.  For any $i\in I$, and using $f,g\in N(\W,\ze)$, we have $Z_i(fg)=(Z_if)g\sub W_ig=\{z_i\}$; thus, $Z_i(fg)=\{z_i\}$, and so $fg=\ze$.
\epf

Now that we have characterised the null subsemigroups of $\T_X$, we wish to find the maximum size of such a subsemigroup.  We begin by calculating the size of the semigroups $N(\W,\ze)$.  It follows from the next lemma that this size depends only on the partitions ${\Z=\set{Z_i}{i\in I}}$ and $\W=\set{W_i}{i\in I}$, indeed only on the cardinalities $|W_i|$ and $|Z_i\sm W_i|$ for each $i\in I$.

\begin{lemma}\label{lem:SizeNW}
If $\W=\set{W_i}{i\in I}$ is a $\ze$-system, then
\[
|N(\W,\ze)| = \prod_{i\in I}|W_i|^{|Z_i\sm W_i|}.
\]
\end{lemma}

\pf
An element $f$ of $N(\W,\ze)$ is determined by the restrictions $f|_{Z_i\sm W_i}$ ($i\in I$), which are arbitrary functions $Z_i\sm W_i\to W_i$.
\epf

If $X$ is infinite, then $\T_X$ contains a null subsemigroup of maximum conceivable size ${2^{|X|}=|\T_X|}$; indeed it follows from Lemma \ref{lem:SizeNW} that $|N(\W,\ze)|=2^{|X|}$ if for example $|W_i|=|Z_i\sm W_i|=|X|$ for some $i$.

We now consider the finite case.  For a positive integer $n$, we write $\nu(n)$ for the maximum size of a null subsemigroup of $\T_n$.  We also write $\nu_r(n)$ for the maximum size of a null subsemigroup of~$\T_n$ whose zero element has rank $r$.

Consider an idempotent $\ze=\trans{Z_1&\cdots&Z_r\\z_1&\cdots&z_r}$ from $\T_n$, and a $\ze$-system $\{W_1,\ldots,W_r\}$ with each $z_i\in W_i$.  We then have two compositions $\si=(s_1,\ldots,s_r)$ and $\tau=(t_1,\ldots,t_r)$, where $s_i=|Z_i|$ and $t_i=|W_i|$ for each~$i$, and we note that $\si\vDash n$ and $\tau\pre\si$, in the notation of Subsection \ref{subsect:comp}.  Furthermore, Lemma \ref{lem:SizeNW} tells us that $|N(\W,\ze)|=\eta(\si,\tau)$.  It follows from the definitions, and using Lemma~\ref{lem:Xi}, that
\begin{align*}
\nu(n) &= \max\set{\eta(\si,\tau)}{\si\vDash n,\ \tau\pre\si}  &\text{and}&& \nu_r(n) &= \max\set{\eta(\si,\tau)}{\si\vDash n,\ |\si|=r,\ \tau\pre\si}\\
&= \max\set{\xi(\si)}{\si\vDash n}  &&&& =\max\set{\xi(\si)}{\si\vDash n,\ |\si|=r} \\
&= \Xi(n) = \xi(n) &&&& = \Xi_r(n) = \xi(n-r+1).
\end{align*}
Combining this with Lemma \ref{lem:Xi}\ref{Xi3}, we have proved the following.

\begin{thm}\label{thm:null}
Let $n$ be a positive integer.
\ben
\item \label{null1} For any $1\leq r\leq n$, the maximum size of a null subsemigroup of $\T_n$ whose zero element has rank $r$ is equal to $\nu_r(n)=\xi(n-r+1)$.  
\item \label{null2}  We have $\nu_1(n) > \nu_2(n) >\cdots>\nu_{n-1}(n)=\nu_n(n) = 1$.
\item \label{null3}  The maximum size of a null subsemigroup of $\T_n$ is equal to $\nu(n) = \nu_1(n) = \xi(n)$.  Moreover, writing $t=\al(n)$, the set
\[
N_n = \bigset{f\in\T_n}{\im(f)\sub\bt \text{ and } \bt f=\{1\}}
\]
is a null subsemigroup of $\T_n$ of size $\nu(n)$.  \epfres
\een
\end{thm}

\begin{rem}
Although the numbers $\nu_r(n)$ are defined in terms of the two parameters $n$ and $r$, there really is only one independent parameter, namely the difference $n-r$.
(Thus, for example, $\nu_5(10) = \nu_{105}(110)=\xi(6)=27$.)  This can be understood in the following way.  Recall from the proof of Lemma \ref{lem:Xi} that the maximum value of $\xi(\si)$, for $\si\vDash n$ with ${|\si|=r}$, occurs when $\si=(n-r+1,1,\ldots,1)$.  This means that the biggest null subsemigroup whose zero~$\ze$ has rank~$r$ occurs when there is a single $\ker(\ze)$-class of size $n-r+1$, 
say $\ze=\trans{\{1,\ldots,n-r+1\}&n-r+2&\cdots&n\\1&n-r+2&\cdots&n}$.  A null subsemigroup $N$ of $\T_n$ with zero $\ze$, and of maximum size $|N|=\nu_r(n)=\xi(n-r+1)$, is then given by $N_{n-r+1}\op\T_{\{n-r+2\}}\op\cdots\op\T_{\{n\}}$, where the subsemigroups $N_k\sub\T_k$ are as defined in Theorem \ref{thm:null}\ref{null3}.  And of course $N$ is then isomorphic to $N_{n-r+1}$.
\end{rem}

\begin{rem}
Null semigroups played an important role in the investigation of subsemigroup chains in $\T_n$ undertaken in \cite{CGMP2017}.  The null subsemigroups in \cite{CGMP2017} were contained in \emph{principal factors} of $\D$-classes, so in a sense all their elements had the same rank (apart from the adjoined zero element).  By contrast, the elements of the null subsemigroups of $\T_n$ considered here can have a range of ranks. 
\end{rem}

Turning Theorem \ref{thm:null}\ref{null3} around, we obtain the degree of a null semigroup:

\begin{thm}\label{thm:null2}
If $S$ is a null semigroup of size $p$, then $\mu(S) = \min\set{n}{\xi(n)\geq p}$.  \epfres
\end{thm}

\begin{rem}\label{rem:LZnull}
Comparing Theorems \ref{thm:null2} and \ref{thm:LZRZ}\ref{LZRZ1}, we see that null and left zero semigroups of the same size have the same degree.
\end{rem}

\subsection{\boldmath Maximal null subsemigroups of $\T_n$}\label{subsect:maxnull}

It is also possible to calculate the number of maximal null subsemigroups of $\T_n$.  (By a maximal null subsemigroup we mean a null subsemigroup that is not properly contained in another null subsemigroup.)  To this end, we write $F(n)$ for this number.  

By Lemma \ref{lem:NWz}, any maximal null subsemigroup of $\T_n$ is of the form $N(\W,\ze)$ for some idempotent~$\ze$, and some $\ze$-system $\W$.  So a first attempt to calculate $F(n)$ would be to count the number of pairs $(\ze,\W)$ where $\ze$ is an idempotent and $\W$ a $\ze$-system.  However, due to a technicality explained below, this results in some over-counting.  To get around this, it will be convenient to also write $F_1(n)$ for the number of maximal null subsemigroups of $\T_n$ whose zero has rank $1$.

Consider an idempotent $\ze=\trans{Z_1&\cdots&Z_r\\z_1&\cdots&z_r}\in\T_n$, and a $\ze$-system $\W=\{W_1,\ldots,W_r\}$, with each $z_i\in W_i\sub Z_i$.  For each $i$, let $\ze_i=\trans{Z_i\\z_i}\in\T_{Z_i}$ be the constant map with image $z_i$.  Then each $\W_i=\{W_i\}$ is trivially a $\ze_i$-system, and it is clear that
\[
N(\W,\ze) = \bigoplus_{i\in\br}N(\W_i,\ze_i),
\]
with each factor $N(\W_i,\ze_i)$ a null subsemigroup of $\T_{Z_i}$ whose zero has rank $1$.  Thus, we can uniquely specify a maximal null subsemigroup $S$ of $\T_n$ as follows:
\ben
\item \label{n1} Choose an integer partition $\si=(s_1,\ldots,s_r)\vdash n$ (which, recall, means that $\si$ is a composition of $n$ and $s_1\geq\cdots\geq s_r$).
\item \label{n2} Choose a set partition $\{Z_1,\ldots,Z_r\}$ of $\bn$ with $|Z_i|=s_i$ for each $i$.
\item \label{n3} Choose maximal null subsemigroups $S_i$ of $\T_{Z_i}$, $i\in\br$, each of which has a zero of rank $1$.
\een
The semigroup so-determined is then $S=S_1\op\cdots\op S_r$.

Given a partition $\si$ as in \ref{n1}, write $f(\si)$ for the number of ways to choose the set partition in~\ref{n2}, and $g(\si)$ for the number of ways to choose the semigroups $S_i$ in \ref{n3}.  Note that the number of choices available in \ref{n3} does not depend on the particular choice made in \ref{n2}.  It is therefore clear that
\begin{equation}\label{eq:F1}
F(n) = \sum_{\si\vdash n}f(\si)g(\si),
\end{equation}
and also that
\begin{equation}\label{eq:F2}
g(\si) = F_1(s_1)\cdots F_1(s_r) \qquad \text{for $\si=(s_1,\ldots,s_r)\vdash n$.}
\end{equation}
A formula for $f(\si)$ is also well known, and easily calculated.  For each $i\in\bn$, let $\si_i$ be the number of entries of $\si$ equal to $i$.  Then
\begin{equation}\label{eq:F3}
f(\si) = \frac{n!}{\prod_{i\in\bn}\si_i!(i!)^{\si_i}}.
\end{equation}
Thus, given \eqref{eq:F1}--\eqref{eq:F3}, to compute $F(n)$ it remains to understand the numbers $F_1(n)$, and we now show that
\begin{equation}\label{eq:F4}
F_1(n) = \begin{cases}
1 &\text{if $n=1$}\\
2 &\text{if $n=2$}\\
n(2^{n-1}-2) &\text{if $n\geq3$.}
\end{cases}
\end{equation}
The $n\leq2$ cases are easily dealt with, so we now assume that $n\geq3$.  
Clearly $F_1(n)$ is equal to~$n$ times the number of maximal null subsemigroups of $\T_n$ whose zero is $\ze=\trans{\bn\\1}$, the constant map with image $\{1\}$, so it is enough to show that there are $2^{n-1}-2$ such subsemigroups.  As above, such a subsemigroup has the form $N(\W,\ze)$, where $\W=\{W\}$ for some subset $\{1\}\sub W\sub\bn$.  For simplicity we denote this subsemigroup by
\[
N(W) = \bigset{f\in\T_n}{\bn f\sub W,\ Wf=\{1\}}.
\]
As extreme cases, we observe that $N(\bn) = N(\{1\}) = \{\ze\}$ is not maximal, as it is obviously contained in every other $N(W)$; recall that we are assuming $n\geq3$.  On the other hand, given distinct subsets ${\{1\}\subset W_1,W_2\subset \bn}$, the subsemigroups $N(W_1)$ and $N(W_2)$ are incomparable in the inclusion order.  Indeed, if ${x\in W_1\sm W_2}$, then for any $y\in\bn\sm W_1$ and $z\in W_2\sm\{1\}$ we have
\begin{equation}\label{eq:W1W2}
\trans{y\ &\bn\sm\{y\}\\x\ &1} \in N(W_1)\sm N(W_2) \AND \trans{x\ &\bn\sm\{x\}\\z\ &1} \in N(W_2)\sm N(W_1).
\end{equation}
It follows from all of this that the maximal null subsemigroups of $\T_n$ with zero $\ze=\trans{\bn\\1}$ are precisely those of the form $N(W)$ with $\{1\} \subset W\subset \bn$.  Since there are of course $2^{n-1}-2$ of these, this completes the proof of \eqref{eq:F4}.

Calculating the numbers $F(n)$ using \eqref{eq:F1}--\eqref{eq:F4} requires the enumeration of all integer partitions of $n$, which is computationally challenging.  However, $F(n)$ may be computed more efficiently using a simple recurrence, as we now explain.  Observe that we can also specify a maximal null subsemigroup $S$ of $\T_n$ in the following way:
\bit
\item Choose the kernel-class containing $1$ of the zero of $S$.  Denote this class by $Z_1$.
\item Choose a maximal null subsemigroup $S_1$ of $\T_{Z_1}$ whose zero has rank $1$.
\item Choose a maximal null subsemigroup $S_2$ of $\T_{\bn\sm Z_1}$ (with zero of any rank).
\eit
The semigroup so-determined is then $S=S_1\op S_2$.  If $|Z_1|=k$, then we can perform the above steps in $\binom{n-1}{k-1}$, $F_1(k)$ and $F(n-k)$ ways, respectively.  This then leads to the recurrence:
\begin{equation}\label{eq:F5}
F(0) = 1 \AND F(n) = \sum_{k=1}^n \binom{n-1}{k-1}F_1(k)F(n-k) \quad\text{for $n\geq1$,}
\end{equation}
with the numbers $F_1(k)$ given in \eqref{eq:F4}.
Table \ref{tab:F} gives some values of $F(n)$ and $F_1(n)$, computed using~\eqref{eq:F4} and \eqref{eq:F5}.  The numbers $F(n)$ grow rapidly; for example, $F(1000)$ has 2287 digits.  At the time of writing, the numbers $F(n)$ did not appear on the OEIS \cite{OEIS}.  The numbers $F_1(n)$ are \cite[Sequence A052749]{OEIS}.

\begin{table}[ht]
\begin{center}
\begin{tabular}{|c|ccccccccccc|}
\hline
$n$ & 0 & 1 & 2 & 3 & 4 & 5 & 6 & 7 & 8 & 9 & 10 \\
\hline
$F(n)$ & 1 & 1 & 3 & 13 & 73 & 451 & 3211 & 26097 & 236433 & 2335123 & 24943171\\
\hline
$F_1(n)$ &  & 1 & 2 & 6 & 24 & 70 & 180 & 434 & 1008 & 2286 & 5100\\
\hline
\end{tabular}
\caption{Values of $F(n)$ and $F_1(n)$:  $F(n)$ is the number of maximal null subsemigroups of $\T_n$, while~$F_1(n)$ is the number of such subsemigroups whose zero has rank $1$.}
\label{tab:F}
\end{center}
\end{table}

\subsection{Nilpotent semigroups}\label{subsect:nil}

The construction of the null semigroups $N(\W,\ze)$ in Subsection \ref{subsect:null} can be generalised to nilpotent semigroups.  Recall that a semigroup $S$ is \emph{nilpotent} if for some integer $k\geq1$ we have $S^k=\{z\}$ for some fixed element $z\in S$, which is necessarily a zero element of $S$; if $k$ is minimal with respect to this property, we say that $S$ is \emph{$k$-nilpotent}.  So $2$-nilpotent semigroups are null.  Here we give the details for $3$-nilpotent subsemigroups of $\T_X$; the construction for $k\geq4$ is easily adapted, though the enumeration becomes more complicated.

To describe $3$-nilpotent subsemigroups of $\T_X$ with zero $\ze=\binom{Z_i}{z_i}_{i\in I}$, we require two systems ${\V=\set{V_i}{i\in I}}$ and $\W=\set{W_i}{i\in I}$ such that $z_i\in V_i\sub W_i\sub Z_i$ for all $i$.  For such a pair of systems, it is easy to show that
\[
N(\V,\W,\ze) = \bigset{f\in\T_X}{Z_if\sub W_i,\ W_if\sub V_i,\ V_if=\{z_i\} \text{ for all $i\in I$}}
\]
is a $3$-nilpotent subsemigroup of $\T_X$ with zero $\ze$.  Strictly speaking, to ensure that $N(\V,\W,\ze)$ is indeed $3$-nilpotent, and not degenerately $2$-nilpotent (i.e., null), we must have \emph{strict} inclusions $\{z_i\}\subset V_i\subset W_i\subset Z_i$ for at least one $i\in I$.  In particular, this forces $|X|\geq4$.
Conversely, any $3$-nilpotent subsemigroup $S$ with zero~$\ze$ is contained in some such $N(\V,\W,\ze)$.  Indeed, given $S$ we define 
\[
W_i = \bigcup_{f\in S}Z_if \AND V_i = \bigcup_{f\in S}W_if \qquad\text{for each $i\in I$.}
\]

The size of $N(\V,\W,\ze)$, as defined above, is given by
\[
|N(\V,\W,\ze)| = \prod_{i\in I}|W_i|^{|Z_i\sm W_i|}|V_i|^{|W_i\sm V_i|}.
\]
In the finite case, say $|X|=n$, writing $\ze=\trans{Z_1&\cdots&Z_r\\z_1&\cdots&z_r}$, $|Z_i|=s_i$, $|W_i|=t_i$ and $|V_i|=u_i$, we have three compositions $\upsilon=(u_1,\ldots,u_r)$, $\tau=(t_1,\ldots,t_r)$ and $\si=(s_1,\ldots,s_r)$, satisfying $\upsilon\pre\tau\pre\si\vDash n$.  We then calculate
\[
|N(\V,\W,\ze)| = t_1^{s_1-t_1}\cdots t_r^{s_r-t_r} \cdot u_1^{t_1-u_1} \cdots u_r^{t_r-u_r} = \eta(\si,\tau)\cdot\eta(\tau,\upsilon).
\]
Now write $u=u_1+\cdots+u_r$ and $t=t_1+\cdots+t_r$, and note that $n=s_1+\cdots+s_r$.  For integers $1\leq i\leq p$ and $1\leq j\leq q$ we have $i^{p-i}j^{q-j}\leq (i+j)^{p-i}(i+j)^{q-j}=(i+j)^{(p+q)-(i+j)}$.  It quickly follows that
\[
|N(\V,\W,\ze)| \leq t^{n-t}\cdot u^{t-u}.
\]
Note that when $r=1$ we have size $t^{n-t}\cdot u^{t-u}$, so the maximum is achieved at rank 1: i.e., a $3$-nilpotent subsemigroup of $\T_n$ of maximum size has a zero of rank $1$.  Thus, if we write~$\ka(n)$ for the maximum size of a $3$-nilpotent subsemigroup of $\T_n$, then we have
\[
\ka(n) = \max\set{t^{n-t}\cdot u^{t-u}}{1\leq u\leq t\leq n}.
\]
For fixed $1\leq t\leq n$, note that $t^{n-t}\cdot u^{t-u}$ is maximised (for $1\leq u\leq t$) when $u=\al(t)$, and then $u^{t-u} = \xi(t)$; the functions $\al$ and $\xi$ were defined in \eqref{eq:xial}.  We have therefore proved the following:

\begin{thm}\label{thm:3nil}
For an integer $n\geq4$, the maximum size of a $3$-nilpotent subsemigroup of $\T_n$ is equal to
\[
\epfreseq
\ka(n) = \max\set{t^{n-t}\cdot\xi(t)}{t\in\bn}.
\]
\end{thm}

Table \ref{tab:3nil} gives calculated values of $\ka(n)$, along with the values of $u$ and $t$ for which ${\ka(n)=t^{n-t}\cdot u^{t-u}}$.  These numbers do not appear on the OEIS \cite{OEIS}.  Figure \ref{fig:nil} shows the egg-box diagram of a $3$-nilpotent subsemigroup of $\T_5$ of size $\ka(5)=18$, again produced with GAP \cite{GAP,Semigroups}.

\begin{table}[ht]
\begin{center}
\scalebox{0.9}{
\begin{tabular}{|c|c|c|r|}
\hline
$n$ & $u$ & $t$ & \multicolumn{1}{c|}{$\ka(n)$} \\
\hline
1 &  &  &  \\
2 &  &  &  \\
3 &  &  &  \\
4 & 2 & 3 & 6 \\
5 & 2 & 3 & 18 \\
6 & 2 & 4 & 64 \\
7 & 2 & 4 & 256 \\
8 & 3 & 5 & 1125 \\
9 & 3 & 6 & 5832 \\
10 & 3 & 6 & 34992 \\
\hline
11 & 3 & 6 & 209952 \\
12 & 3 & 7 & 1361367 \\
13 & 3 & 7 & 9529569 \\
14 & 4 & 8 & 67108864 \\
15 & 4 & 9 & 544195584 \\
16 & 4 & 9 & 4897760256 \\
17 & 4 & 9 & 44079842304 \\
18 & 4 & 10 & 409600000000 \\
19 & 4 & 10 & 4096000000000 \\
20 & 4 & 10 & 40960000000000 \\
\hline
\end{tabular}
\qquad\qquad
\begin{tabular}{|c|r|r|r|}
\hline
$n$ & $u$ & $t$ & \multicolumn{1}{c|}{$\ka(n)$} \\
\hline
21 & 4 & 11 & 424958764662784 \\
22 & 5 & 12 & 4837294080000000 \\
23 & 5 & 12 & 58047528960000000 \\
24 & 5 & 13 & 700062653920703125 \\
25 & 5 & 13 & 9100814500969140625 \\
26 & 5 & 13 & 118310588512598828125 \\
27 & 5 & 14 & 1550224166512000000000 \\
28 & 5 & 14 & 21703138331168000000000 \\
29 & 5 & 14 & 303843936636352000000000 \\
30 & 6 & 15 & 4412961507515625000000000 \\
\hline
31 & 6 & 16 & 69712754611742420055883776 \\
32 & 6 & 16 & 1115404073787878720894140416 \\
33 & 6 & 16 & 17846465180606059534306246656 \\
34 & 6 & 17 & 300120331617031984667981862912 \\
35 & 6 & 17 & 5102045637489543739355691669504 \\
36 & 6 & 17 & 86734775837322243569046758381568 \\
37 & 6 & 18 & 1541674189500358697210578156388352 \\
38 & 6 & 18 & 27750135411006456549790406814990336 \\
39 & 7 & 19 & 520293618503860588010973963362962801 \\
40 & 7 & 20 & 10159549097653043200000000000000000000 \\
\hline
\end{tabular}
}
\caption{Values of $\ka(n)$, which is the maximum size of a $3$-nilpotent subsemigroup of $\T_n$.  Also shown are the values of $u$ and $t$ for which $\ka(n)=t^{n-t}\cdot u^{t-u}$.}
\label{tab:3nil}
\end{center}
\end{table}

\begin{figure}[!ht]
\begin{center}
\includegraphics[width=\textwidth]{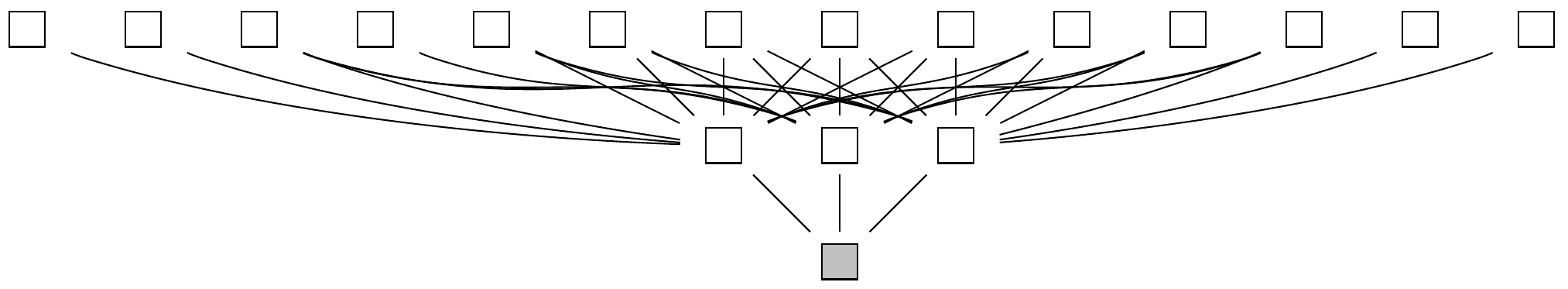}
\caption{Egg-box diagram of a $3$-nilpotent subsemigroup of $\T_5$ of maximum possible size $\ka(5)=18$.}
\label{fig:nil}
\end{center}
\end{figure}

\subsection{\boldmath Maximal $3$-nilpotent subsemigroups of $\T_n$}\label{subsect:maxnil}

As with null semigroups in Subsection \ref{subsect:maxnull}, it is possible to calculate the number of maximal $3$-nilpotent subsemigroups of $\T_n$, though the calculation is a little more delicate.  Here, by a maximal $3$-nilpotent subsemigroup we mean a $3$-nilpotent subsemigroup that is not properly contained in another $3$-nilpotent subsemigroup.  To this end, we write
\bit
\item $G(n)$ for the number of maximal $3$-nilpotent subsemigroups of $\T_n$, and
\item $G_1(n)$ for the number of such subsemigroups whose zero has rank $1$.
\eit
For reasons that will become clear shortly, we also need to define some further sequences.  Specifically, we also write
\bit
\item $H(n)$ for the number of null subsemigroups of $\T_n$ that are not properly contained in any null or $3$-nilpotent subsemigroup, and
\item $H_1(n)$ for the number of such subsemigroups whose zero has rank $1$.
\eit
Now consider some maximal $3$-nilpotent subsemigroup $S$ of $\T_n$.  So this must have the form ${S=N(\V,\W,\ze)}$, as described in the previous subsection, and we write
\[
\ze=\trans{Z_1&\cdots&Z_r\\z_1&\cdots&z_r} \COMMA {\W=\{W_1,\ldots,W_r\}} \AND \V=\{V_1,\ldots,V_r\},
\]
where 
\begin{equation}\label{eq:zVWZ}
\{z_i\}\sub V_i\sub W_i\sub Z_i \qquad\text{for each $i$.}
\end{equation}
Again, the inclusions in \eqref{eq:zVWZ} do not all have to be strict, but to ensure that $S$ is $3$-nilpotent (and not null) there must be at least one $i$ such that $\{z_i\}\subset V_i\subset W_i\subset Z_i$.  For each $i\in\br$, let $\ze_i=\trans{Z_i\\z_i}\in\T_{Z_i}$, $\W_i=\{W_i\}$ and $\V_i=\{V_i\}$, so that
\[
S = \bigoplus_{i\in\br} S_i \qquad\text{where each}\qquad S_i=N(\V_i,\W_i,\ze_i).
\]
Moreover, since $S$ is maximal $3$-nilpotent, each $S_i=N(\V_i,\W_i,\ze_i)$ is either:
\bit
\item a maximal $3$-nilpotent subsemigroup of $\T_{Z_i}$, or else 
\item a null subsemigroup of $\T_{Z_i}$ that is not properly contained in any null or $3$-nilpotent subsemigroup,
\eit
depending on the strictness of the inclusions in \eqref{eq:zVWZ} for that particular $i\in\br$.  In either case, the zero $\ze_i$ of each factor $S_i$ has rank $1$.  One can then follow either of the two methods in Subsection~\ref{subsect:maxnull} to calculate the number $G(n)$ of maximal $3$-nilpotent subsemigroups of $\T_n$.  

We just give the details for the second approach, and obtain a recurrence for $G(n)$.  
The basis of the recursion is the observation that a maximal $3$-nilpotent subsemigroup~$S$ can be specified in the following way:
\ben
\item \label{3n1} Choose the kernel-class containing $1$ of the zero of $S$.  Denote this class by $Z_1$.
\item \label{3n2} Choose a subsemigroup $S_1$ of $\T_{Z_1}$ whose zero has rank $1$, and is either
\bit
\item maximal $3$-nilpotent, or else
\item null, but not properly contained in any null or $3$-nilpotent subsemigroup of $\T_{Z_1}$.
\eit
\item \label{3n3} Choose a subsemigroup $S_2$ of $\T_{\bn\sm Z_1}$ (with zero of any rank) that is either
\bit
\item maximal $3$-nilpotent, or else
\item null, but not properly contained in any null or $3$-nilpotent subsemigroup of $\T_{\bn\sm Z_1}$,
\eit
ensuring that $S_1$ and $S_2$ are not both null.
\een
The semigroup so-determined is then $S=S_1\op S_2$.  If $|Z_1|=k$, then we can perform step \ref{3n1} in~$\binom{n-1}{k-1}$ ways.  If $S_1$ is null, then there are $H_1(k)$ and $G(n-k)$ ways to perform steps \ref{3n2} and~\ref{3n3}, respectively.  Otherwise, $S_1$ is $3$-nilpotent, and there are $G_1(k)$ and $G(n-k)+H(n-k)$ ways.  This then leads to the recurrence:
\begin{align}
\label{eq:G1} G(0)&=0,\\
\label{eq:G2} G(n) &= \sum_{k=1}^n \binom{n-1}{k-1} \Big( H_1(k)G(n-k) + G_1(k)\big(G(n-k)+H(n-k)\big) \Big) \qquad\text{for $n\geq1$.}
\end{align}
In order to implement this, it remains to understand the numbers $G_1(n)$, $H_1(n)$ and $H(n)$.  Formulae/recursions for these are given below in \eqref{eq:H1}--\eqref{eq:G3}.

We begin with the $H$-numbers, and we first claim that
\begin{equation}\label{eq:H1}
H_1(n) = \begin{cases}
n! &\text{if $1\leq n\leq3$}\\
0 &\text{if $n\geq4$.}
\end{cases}
\end{equation}
Indeed, for $n\leq3$ we have $H_1(n)=F_1(n)$, as there are no (properly) $3$-nilpotent subsemigroups of $\T_n$ for $n\leq3$.  We now assume that $n\geq4$.  Aiming for a contradiction, suppose $S$ is a null subsemigroup of $\T_n$, with zero of rank $1$, that is not properly contained in any null or $3$-nilpotent subsemigroup.  By the former, $S$ is a maximal null subsemigroup.  By symmetry, we may assume the zero of $S$ is $\ze=\trans{\bn\\1}$, so that $S=N(\W,\ze)$, where $\W=\{W\}$ for some $\{1\}\subset W\subset\bn$.  Since $n\geq4$, we either have $|W|\geq3$ or else $|W|\leq n-2$.  But in these cases, respectively, we have:
\bit
\item $S=N(\W,\ze)\subset N(\W',\W,\ze)$, where $\W'=\big\{W\sm\{x\}\big\}$ for any $x\in W\sm\{1\}$.
\item $S=N(\W,\ze)\subset N(\W,\W'',\ze)$, where $\W''=\big\{W\cup\{y\}\big\}$ for any $y\in\bn\sm W$.
\eit
Since $N(\W',\W,\ze)$ and $N(\W,\W'',\ze)$ are both $3$-nilpotent, we have reached the desired contradiction, and this completes the proof of \eqref{eq:H1}.

As in Subsection \ref{subsect:maxnull}, it is then easy to see that
\begin{equation}\label{eq:H}
H(0) = 1 \AND H(n) = \sum_{k=1}^n \binom{n-1}{k-1}H_1(k)H(n-k) \quad\text{for $n\geq1$.}
\end{equation}
Keeping \eqref{eq:H1} in mind, we note that only the $k\leq3$ terms in the above sum are non-zero.

As for the $G$-numbers, we claim that
\begin{equation}\label{eq:G3}
G_1(n) = \begin{cases}
0 &\text{if $1\leq n\leq3$}\\
3n(1+3^{n-2}-2^{n-1}) &\text{if $n\geq4$.}
\end{cases}
\end{equation}
As usual, this is clear for $n\leq3$, so we assume that $n\geq4$.  And again, $G_1(n)$ is equal to $n$ times the number of maximal $3$-nilpotent subsemigroups of $\T_n$ whose zero is $\ze=\trans{\bn\\1}$.  So consider some such subsemigroup~$S$, which has the form 
\begin{equation}\label{eq:NVWz}
S = \bigset{f\in\T_n}{\bn f\sub W,\ Wf\sub V,\ Vf=\{1\}} \qquad\text{for some $\{1\}\subset V\subset W\subset \bn$.}
\end{equation}
Here, $S=N(\V,\W,\ze)$, for $\W=\{W\}$ and~$\V=\{V\}$, and we abbreviate this to $S=N(V,W)$.  This time, ${S_1=N(V_1,W_1)}$ and ${S_2=N(V_2,W_2)}$ are incomparable if $(V_1,W_1)\not=(V_2,W_2)$.  (If $V_1\not=V_2$ or $W_1\not=W_2$, then as in \eqref{eq:W1W2}, it is easy to construct transformations belonging to $S_1\sm S_2$ and $S_2\sm S_1$.)  Thus, $G_1(n)/n$ is equal to the number of pairs $(V,W)$, as in \eqref{eq:NVWz}.  Specifying such a pair is equivalent to specifying the ordered partition $\big(V\sm\{1\}, W\sm V, \bn\sm W)$ of $\bn\sm\{1\}$.  It follows that $G_1(n)/n = 3!\cdot S(n-1,3)$, where $S(k,l)$ is a Stirling number of the second kind.  It is well known (and easy to verify using the standard Stirling recurrence) that $S(k,3) = (1+3^{k-1}-2^k)/2$ for $k\geq3$, so \eqref{eq:G3} quickly follows.

Table \ref{tab:G} gives some values of $G(n)$ and $G_1(n)$, computed using \eqref{eq:G1}--\eqref{eq:G3}.  The numbers~$G_1(n)$ appear as Sequence A052761 on the OEIS \cite{OEIS}, while $G(n)$ was not on the OEIS at the time of writing.  Although $F(n)>G(n)$ for $n\leq6$, the sequence $G(n)$ quickly starts to out-grow~$F(n)$.  For example, $G(1000)$ has 2405 digits, and we mentioned earlier than $F(1000)$ has 2287.

\begin{table}[ht]
\begin{center}
\begin{tabular}{|c|ccccccccccc|}
\hline
$n$ & 0 & 1 & 2 & 3 & 4 & 5 & 6 & 7 & 8 & 9 & 10 \\
\hline
$G(n)$ & 0 & 0 & 0 & 0 & 24 & 300 & 3060 & 32340 & 353808 & 4199580 & 54149820\\
\hline
$G_1(n)$ & & 0 & 0 & 0 & 24 & 180 & 900 & 3780 & 14448 & 52164 & 181500\\
\hline
\end{tabular}
\caption{Values of $G(n)$ and $G_1(n)$:  $G(n)$ is the number of maximal $3$-nilpotent subsemigroups of~$\T_n$, while~$G_1(n)$ is the number of such subsemigroups whose zero has rank $1$.}
\label{tab:G}
\end{center}
\end{table}

\subsection{Right null semigroups}\label{subsect:RN}

When studying variants in Section \ref{sect:Tna}, an important role will be played by certain special kinds of semigroups we call \emph{right null semigroups}.

Recall that a \emph{band} is a semigroup consisting entirely of idempotents.  Following Clifford \cite{Clifford1954}, a semigroup $S$ is a \emph{band of semigroups} if we have $S=\bigcup_{b\in B}S_b$, for some band~$B$, where the $S_b$ are pairwise-disjoint subsemigroups of $S$, and $S_bS_c\sub S_{bc}$ for all $b,c\in B$.  
Of particular interest for us is the case in which:
\bit
\item $B$ is a right zero semigroup ($bc=c$ for all $b,c\in B$), 
\item each $S_b$ is a null semigroup, say with zero $z_b\in S_b$, and 
\item $S_bS_c=\{z_{c}\}$ for all $b,c\in B$.
\eit
We will call such a semigroup $S=\bigcup_{b\in B}S_b$ a \emph{right null semigroup}.  When $|B|=1$, $S$ is a null semigroup; when each $|S_b|=1$, $S$ is a right zero semigroup; so right null semigroups simultaneously generalise both classes.

For a right null semigroup $S$ as above, the idempotents $E(S)={\set{z_b}{b\in B}}$ form a right zero subsemigroup of $S$ isomorphic to $B$; so we will typically identify $E(S)$ with $B$.  We call~$S$ \emph{uniform} if the $S_b$ have a common size.  If $|B|=p$ and $|S_b|=q$ for all $b\in B$, we call $S$ a $p\times q$ (uniform) right null semigroup, and then of course $|S|=pq$.  In fact, it is easy to show that a $p\times q$ right null semigroup is simply a direct product of a $p$-element right zero semigroup and a $q$-element null semigroup.

In order to understand the right null subsemigroups of a full transformation semigroup $\T_X$, we make use of the following natural extension of the $\ze$-systems, and associated null semigroups, from Definition \ref{defn:ze}.

\begin{defn}\label{defn:NWB}
Given a right zero subsemigroup $B$ of $\T_X$, a \emph{$B$-system} is a collection $\W$ of sets such that $\W$ is simultaneously a $\ze$-system for all $\ze\in B$.  For such a $B$-system $\W$, we define
\[
N(\W,B) = \bigcup_{\ze\in B}N(\W,\ze).
\]
\end{defn}

\begin{lemma}\label{lem:NWB}
Given a right zero subsemigroup $B$ of $\T_X$, a subset $S$ of $\T_X$ is a right null semigroup with $E(S)=B$ if and only if ${B\sub S\sub N(\W,B)}$ for some $B$-system $\W$.
\end{lemma}

\pf
The proof is easily adapted from that of Lemma \ref{lem:NWz}.  For the forwards implication we still define $W_i=\bigcup_{f\in S}Z_if$.
\epf

Of course we have $|N(\W,B)|=|B|\cdot|N(\W,\ze)|$ for any $\ze\in B$.  Moreover, given a right zero subsemigroup $B\sub\T_X$, one could find the maximum size of a right null subsemigroup with idempotents~$B$, by choosing $\W$ so that each $N(\W,\ze)$, $\ze\in B$, is as large as possible.

On the other hand, given a partition $\Z=\set{Z_i}{i\in I}$ of $X$, and a system $\W=\set{W_i}{i\in I}$ with $\emptyset\not= W_i\sub Z_i$ for all $i\in I$, we can construct a maximum band $B=B(\Z,\W)$ such that $\W$ is a $B$-system.  Specifically, we take $B$ to consist of all idempotents of the form $\ze=\binom{Z_i}{w_i}_{i\in I}$ for all choices of $w_i\in W_i$ ($i\in I$).  This $B$ has size $\prod_{i\in I}|W_i|$.

When $X$ is infinite, it is easy to describe $2^{|X|}\times2^{|X|}$ uniform right null subsemigroups of $\T_X$.

The next result gives necessary and sufficient conditions for a full transformation semigroup~$\T_m$ to contain a right null semigroup of specified dimensions.  The statement uses compositions and their associated parameters, as defined in Subsection \ref{subsect:comp}.

\begin{lemma}\label{lem:mpq}
For $m,p,q\in\N$, $\T_m$ contains a $p\times q$ uniform right null subsemigroup if and only if there exist compositions $\si$ and $\tau$ such that
\begin{equation}\label{eq:mpq}
\tau\pre\si\vDash m \COMMA \textstyle{\prod\tau\geq p} \COMMA \eta(\si,\tau)\geq q.
\end{equation}
\end{lemma}

\pf
By Lemma \ref{lem:NWB}, $\T_m$ contains such a subsemigroup if and only it contains a subsemigroup of the form $N(\W,B)$ for 
\bit
\item some right zero subsemigroup $B\sub\T_m$ with $|B|\geq p$, and 
\item some $B$-system $\W$, such that $|N(\W,\ze)|\geq q$ for all $\ze\in B$.  
\eit
Consider some such $B$ and $\W$.  Suppose the common kernel-classes of the elements of $B$ are $Z_1,\ldots,Z_r$, and write $\W=\{W_1,\ldots,W_r\}$, where $W_i\sub Z_i$ for each $i$.    Also, write $s_i=|Z_i|$ and $t_i=|W_i|$ for all $i$.  Clearly the compositions $\si=(s_1,\ldots,s_r)$ and $\tau=(t_1,\ldots,t_r)$ satisfy $\tau\pre\si\vDash m$.  Moreover, since every $\ze\in B$ has the form $\ze = \trans{Z_1&\cdots&Z_r\\w_1&\cdots&w_r}$ for some $w_i\in W_i$, we have
\[
p \leq |B| \leq t_1\cdots t_r = \textstyle{\prod\tau} \AND q \leq |N(\W,\ze)| = t_1^{s_1-t_1}\cdots t_r^{s_r-t_r} = \eta(\si,\tau).
\]

Conversely, suppose the compositions $\si=(s_1,\ldots,s_r)$ and $\tau=(t_1,\ldots,t_r)$ satisfy conditions~\eqref{eq:mpq}.  We construct a (uniform) $\prod\tau\times\eta(\si,\tau)$ right null subsemigroup $S$ of $\T_m$ as follows, and then~$S$ then of course contains a $p\times q$ right null subsemigroup.  Choose a partition ${\Z=\{Z_1,\ldots,Z_r\}}$ of $\bm$, and a system $\W=\{W_1,\ldots,W_r\}$ such that $W_i\sub Z_i$, $|Z_i|=s_i$ and $|W_i|=t_i$ for each $i$.  We define $S=N(\W,B)$, as in Definition \ref{defn:NWB}, for the band $B=B(\Z,\W)$ described after the proof of Lemma \ref{lem:NWB}.
\epf

It will transpire in Section \ref{sect:Tna} that we need to know the degree of a certain special kind of right null semigroup.  Since all $p\times q$ right null semigroups are isomorphic, we will write $\rn(p,q)$ for the degree of any such semigroup.  

\begin{thm}\label{thm:RN}
\ben
\item \label{RN1} For integers $p,q\geq1$ we have
\[
\rn(p,1) = \min\set{n}{\pi(n)\geq p} \AND \rn(1,q) = \min\set{n}{\xi(n)\geq q}.
\]
\item \label{RN2} For integers $p,q\geq2$ we have 
\begin{align*}
\rn(p,q) &= \min \set{2A+3B+C+D}{0\leq A\leq2,\ B,D\geq0,\ C\geq2,\ 2^A3^BC\geq p,\ C^D\geq q}\\
&= {\min} \bigset{2A + 3\lceil\log_3(p/2^AC)\rceil + C + \lceil\log_Cq\rceil}{ 0\leq A\leq 2\leq C\leq \rn(p,1) + \rn(1,q) },
\end{align*}
where the parameters $A,B,C,D$ in the above sets are all integers.
\een
\end{thm}

\pf 
\firstpfitem{\ref{RN1}} This follows from Theorems \ref{thm:LZRZ}\ref{LZRZ2} and \ref{thm:null2}, given that $p\times1$ and $1\times q$ right null semigroups are right zero or null semigroups, respectively.

\pfitem{\ref{RN2}}  We now assume that $p,q\geq2$.  By Lemma \ref{lem:mpq}, the number $\rn(p,q)$ is the solution to the following optimisation problem:

\begin{oprob}\label{oprob1}
Minimise $s_1+\cdots+s_r$, for all $r,s_1,\ldots,s_r,t_1,\ldots,t_r\in\N$, subject to:
\[
t_i\leq s_i \text{ for all $i$} \COMMA
t_1\cdots t_r\geq p \COMMA
t_1^{s_1-t_1}\cdots t_r^{s_r-t_r}\geq q.
\]
\end{oprob}

In the remainder of the proof, we reformulate Problem \ref{oprob1} until we arrive at the desired result.
Writing $d_i=s_i-t_i$ in the above, Problem \ref{oprob1} is equivalent to the following:

\begin{oprob}\label{oprob2}
Minimise $(t_1+\cdots+t_r)+(d_1+\cdots+d_r)$, for all $r,t_1,\ldots,t_r\in\N$ and ${d_1,\ldots,d_r\in\N\cup\{0\}}$, subject to:
\[
t_1\cdots t_r\geq p \AND t_1^{d_1}\cdots t_r^{d_r}\geq q.
\]
\end{oprob}

In the above, some of the $d_i$ can be $0$, but not all (as $q\geq2$).  For a given sum $d_1+\cdots+d_r=D$, we have $t_1^{d_1}\cdots t_r^{d_r}\leq \max(t_1,\ldots,t_r)^D$, so Problem \ref{oprob2} is equivalent to:

\begin{oprob}\label{oprob3}
Minimise $t_1+\cdots+t_r+D$, for all $r,t_1,\ldots,t_r,D\in\N$, subject to:
\[
t_1\cdots t_r\geq p \AND \max(t_1,\ldots,t_r)^D\geq q.
\]
\end{oprob}

Adding an ordering constraint $t_1\leq\cdots\leq t_r$ of course does not change the solution to Problem~\ref{oprob3}.  Additionally writing $C=t_r=\max(t_1,\ldots,t_r)$, we can simplify the statement a little:

\begin{oprob}\label{oprob4}
Minimise $t_1+\cdots+t_{r-1}+C+D$, for all $r,t_1,\ldots,t_{r-1},C,D\in\N$, subject to:
\[
t_1\leq\cdots\leq t_{r-1}\leq C \COMMA t_1\cdots t_{r-1}C\geq p \COMMA C^D\geq q.
\]
\end{oprob}

Any solution to Problem \ref{oprob4} has $t_1,\ldots,t_{r-1},C\geq2$. Indeed, $C\geq2$ follows from $C^D\geq q\geq2$, while if $r\geq2$ and $t_1=1$, then $t_1$ could be removed from the sum $t_1+\cdots+t_{r-1}+C+D$ without decreasing the product $t_1\cdots t_{r-1}C$.  

Also, for a given sum $t_1+\cdots+t_{r-1}$, the maximum value of $t_1\cdots t_{r-1}$ occurs when at most two of the $t_i$ are $2$, and the rest are $3$.  Problem \ref{oprob4} is then equivalent to:

\begin{oprob}\label{oprob5}
Minimise $2A+3B+C+D$, for integers
\[
0\leq A\leq 2 \COMMA B,D\geq0 \COMMA C\geq2, 
\]
subject to $2^A3^BC\geq p$ and $C^D\geq q$.
\end{oprob}

This gives the first claimed expression for $\rn(p,q)$, and we now work towards establishing the second.  To do so, suppose we have some solution $(A,B,C,D)$ to Problem \ref{oprob5}, so that ${2A+3B+C+D=\mu(S)}$, where $S$ is some $p\times q$ right null semigroup.  Since ${C^D\geq q}$ we must have $D\geq\log_Cq$, and so minimality forces $D=\lceil\log_Cq\rceil$.  We similarly deduce from $2^A3^BC\geq p$ that $B = \lceil\log_3(p/2^AC)\rceil$.  

It remains to show that $C$ is bounded above by ${\rn(p,1)+\rn(1,q)}$.  To see this, note first that $C\leq2A+3B+C+D=\mu(S)$.  But $S\cong R\times N$ for some $p$-element right zero semigroup~$R$ and some $q$-element null semigroup $N$, so $C\leq\mu(S)\leq\mu(R)+\mu(N)=\rn(p,1)+\rn(1,q)$, as required.
\epf

\begin{rem}\label{rem:RN}
When $p,q\geq2$, the second expression in Theorem \ref{thm:RN}\ref{RN2} reduces the calculation of $\rn(p,q)$ to three separate (finite) one-variable minimisations, one for each value of $A=0,1,2$.  It is then easy to compute the numbers $\rn(p,q)$; some values are given in Table \ref{tab:RN}.  We also observe (from the proof) that the upper bound on $C$ of $\rn(p,1)+\rn(1,q)$ could be replaced by any other \emph{a priori} known upper bound on $\rn(p,q)$.

It is interesting to compare Table \ref{tab:RN} with Table \ref{tab:RB}, which gives the degree $\be(p,q)$ of a $p\times q$ rectangular band.  Indeed, consider a $p\times q$ rectangular band $B$ and a $q\times p$ right null semigroup~$S$ (note the swapping of the parameters for $S$).  Then we have
\[
B\cong L\times R \AND S\cong R\times N
\]
for some left zero semigroup $L$, right zero semigroup $R$ and null semigroup $N$, with ${|L|=|N|=p}$ and $|R|=q$.  In fact, not only do we have $|L|=|N|$, but we have already observed in Remark~\ref{rem:LZnull} that $\mu(L)=\mu(N)$.  However, we do not necessarily have $\mu(B)=\mu(S)$; for example, $\be(2,2)=4$ while $\rn(2,2)=3$.  Thus, we have a natural example where
\[
\mu(T\times U)\not=\mu(T\times V) \qquad\text{despite having}\qquad \mu(U)=\mu(V).
\]
Other examples of this phenomenon are not hard to find, and some can be constructed using other results from this paper.  For example, let $B$ be $2\times2$ rectangular band, and let $L$ be a left zero semigroup and $G$ a group with $|L|=|G|=2$.  Clearly $\mu(L)=\mu(G)=2$.  Consulting Table~\ref{tab:RB}, we have $\mu(B)=4=\be_2(2,2)$.  It then follows from Theorem \ref{thm:RG} that $\mu(B\times G)=4$ as well.  On the other hand, $B\times L$ is a $4\times2$ rectangular band, and so $\mu(B\times L)=5$, again by Table \ref{tab:RB}.
\end{rem}

\begin{rem}\label{rem:RN2}
The above results also lead to a formula for the degree of an arbitrary finite right null semigroup.  Specifically, if $S=\bigcup_{b\in B}S_b$ is some such semigroup, then
\[
\mu(S)=\rn(p,q) \qquad\text{where}\qquad p=|B| \text{ and } q={\max}\bigset{|S_b|}{b\in B}.
\]
Indeed, it follows from Lemma \ref{lem:NWB} that if $S$ embeds in $\T_n$ then $S$ is contained in some uniform right null subsemigroup of $\T_n$, which of course contains a $p\times q$ right null semigroup.
\end{rem}

\begin{rem}\label{rem:LN}
Of course one could also consider \emph{left} null semigroups.  These are defined as expected, and the uniform ones are of course direct products of left zero and null semigroups.  A description of left null \emph{transformation} semigroups is more complicated than for right null semigroups, and is omitted since it is not important for our central purposes (cf.~Subsection~\ref{subsect:r=1}).  As an example of the difference in behaviour between left and right null transformation semigroups, the following can all be easily verified using GAP \cite{GAP,Semigroups}:
\bit
\item The transformations $\trans{1&2&3&4&5\\1&1&x&5&5}$, $x=1,2,4,5$, form a $2\times2$ (uniform) left null subsemigroup of $\T_5$.
\item The semigroup $\T_4$ contains no $2\times2$ left null subsemigroup.  
\item The transformations $\trans{1&2&3&4\\1&1&x&4}$, $x=1,2,4$, form a left null subsemigroup of $\T_4$ with two idempotents.
\eit
The first two points show that the degree of a $2\times2$ left null semigroup is $5$.  Combined with the third point, it follows that the degree of an arbitrary left null semigroup can not be reduced to the uniform case, contrasting with the right null situation (cf.~Remark \ref{rem:RN2}).
\end{rem}

\begin{table}[htp]
\begin{center}
\begin{tabular}{|l|cccccccccc|}
\hline
$	p\sm q	$ & $	1	$ & $	2	$ & $	3	$ & $	4	$ & $	5	$ & $	6	$ & $	7	$ & $	8	$ & $	9	$ & $	10	$ \\
\hline
$\phantom{1}1$ & 1 & 3 & 4 & 4 & 5 & 5 & 5 & 5 & 5 & 6 \\ 
$\phantom{1}2$ & 2 & 3 & 4 & 4 & 5 & 5 & 5 & 5 & 5 & 6 \\
$\phantom{1}3$ & 3 & 4 & 4 & 5 & 5 & 5 & 5 & 5 & 5 & 6 \\ 
$\phantom{1}4$ & 4 & 5 & 5 & 5 & 6 & 6 & 6 & 6 & 6 & 6 \\
$\phantom{1}5$ & 5 & 6 & 6 & 6 & 6 & 7 & 7 & 7 & 7 & 7 \\ 
$\phantom{1}6$ & 5 & 6 & 6 & 7 & 7 & 7 & 7 & 7 & 7 & 8 \\
$\phantom{1}7$ & 6 & 7 & 7 & 7 & 8 & 8 & 8 & 8 & 8 & 8 \\ 
$\phantom{1}8$ & 6 & 7 & 7 & 7 & 8 & 8 & 8 & 8 & 8 & 8 \\
$\phantom{1}9$ & 6 & 7 & 7 & 8 & 8 & 8 & 8 & 8 & 8 & 9 \\ 
$10$ & 7 & 8 & 8 & 8 & 8 & 9 & 9 & 9 & 9 & 9 \\
\hline
\end{tabular}
\\[5mm]
\begin{tabular}{|l|cccccccccc|}
\hline
$	\phantom{1}p\sm q	$ & $	10	$ & $	20	$ & $	30	$ & $	40	$ & $	50	$ & $	60	$ & $	70	$ & $	80	$ & $	90	$ & $	100	$ \\
\hline
$\phantom{1}10$ & \phantom{1}9 & \phantom{1}9 & 10 & 10 & 10 & 10 & 10 & 10 & 10 & 10 \\
$\phantom{1}20$ & 11 & 11 & 12 & 12 & 12 & 12 & 12 & 12 & 12 & 12 \\
$\phantom{1}30$ & 12 & 12 & 13 & 13 & 13 & 13 & 13 & 13 & 13 & 13 \\
$\phantom{1}40$ & 13 & 13 & 14 & 14 & 14 & 14 & 14 & 14 & 14 & 14 \\
$\phantom{1}50$ & 14 & 14 & 14 & 15 & 15 & 15 & 15 & 15 & 15 & 15 \\
$\phantom{1}60$ & 14 & 14 & 15 & 15 & 15 & 15 & 15 & 15 & 15 & 15 \\
$\phantom{1}70$ & 14 & 15 & 15 & 15 & 15 & 15 & 16 & 16 & 16 & 16 \\
$\phantom{1}80$ & 15 & 15 & 16 & 16 & 16 & 16 & 16 & 16 & 16 & 16 \\
$\phantom{1}90$ & 15 & 15 & 16 & 16 & 16 & 16 & 16 & 16 & 16 & 16 \\
$100$ & 15 & 16 & 16 & 16 & 16 & 16 & 17 & 17 & 17 & 17 \\
\hline
\end{tabular}
\\[5mm]
\begin{tabular}{|l|cccccccccc|}
\hline
$	p\sm q	$ & $	10^0	$ & $	10^1	$ & $	10^2	$ & $	10^3	$ & $	10^4	$ & $	10^5	$ & $10^6	$ & $	10^7	$ & $	10^8	$ & $10^9$	 \\
\hline
$	10^0	$ & \phantom{1}1 & \phantom{1}6 & \phantom{1}8 & \phantom{1}9 & 11 & 13 & 14 & 15 & 17 & 18 \\
$	10^1	$ & \phantom{1}7 & \phantom{1}9 & 10 & 12 & 13 & 15 & 16 & 17 & 18 & 19 \\
$	10^2	$ & 13 & 15 & 17 & 18 & 20 & 21 & 22 & 23 & 25 & 26 \\
$	10^3	$ & 20 & 22 & 23 & 25 & 26 & 27 & 29 & 30 & 31 & 32 \\
$	10^4	$ & 26 & 28 & 29 & 31 & 32 & 33 & 35 & 36 & 37 & 38 \\
$	10^5	$ & 32 & 34 & 36 & 37 & 39 & 40 & 41 & 42 & 43 & 44 \\
$	10^6	$ & 38 & 41 & 42 & 43 & 45 & 46 & 47 & 48 & 50 & 51 \\
$	10^7	$ & 45 & 47 & 48 & 50 & 51 & 52 & 54 & 55 & 56 & 57 \\
$	10^8	$ & 51 & 53 & 55 & 56 & 57 & 58 & 60 & 61 & 62 & 63 \\
$	10^9	$ & 57 & 59 & 61 & 62 & 64 & 65 & 66 & 67 & 68 & 69 \\
\hline
\end{tabular}
\caption{Calculated values of $\rn(p,q)$, which is the degree of a $p\times q$ right null semigroup.}
\label{tab:RN}
\end{center}
\end{table}

\subsection{Decreasing the degree by forming a variant}\label{subsect:Sa<S}

Before we move on, we pause to address another question asked in \cite{JE2020}.  

Consider a finite semigroup $S$, and an arbitrary element $a\in S$.  It is natural to wonder how the degrees of $S$ and the variant $S^a$ are related.  Theorem \ref{thm:EM} leads to the upper bound ${\mu(S^a)<2\mu(S)}$.  Indeed, if $n=\mu(S)$, then any embedding $\phi:S\to\T_n$ yields an embedding ${S^a\to\T_n^{a\phi}}$, which we can follow with an embedding ${\T_n^{a\phi}\to\T_{2n-r}}$, where $r=\rank(a\phi)$.  

Question~4.4 of \cite{JE2020} asks if there exists a finite semigroup $S$ such that $\mu(S^a)<\mu(S)$ for some $a\in S$.  This can be quickly answered in the affirmative, and we do so with two examples.  The first involves cyclic groups, so we begin with the following basic result.  It is most likely well known, but proofs are provided for convenience.

\begin{lemma}\label{lem:FG}
Let $G$ be a finite group, and let $\phi:G\to\T_n$ be an embedding, where $n=\mu(G)$.
\ben
\item \label{FG1} The image $\im(\phi)$ is contained in the symmetric group $\S_n$.
\item \label{FG2} For every $i\in\bn$ there exists $g\in G$ such that $i(g\phi)\not=i$.
\item \label{FG3} If $S=G\cup\{0\}$ is the semigroup obtained by adjoining a zero element to $G$, then ${\mu(S)=\mu(G)+1}$.
\een
\end{lemma}

\pf
\firstpfitem{\ref{FG1}}  Clearly $\im(\phi)$ is contained in some group $\H$-class $H$ of $\T_n$, and $H$ is contained in some $\D$-class $D_r$ with $1\leq r\leq n$.  But $H\cong\S_r$, so $G$ embeds in $\S_r$, and hence in $\T_r$.  Thus, by minimality of $n=\mu(G)$ we have $n\leq r$, and so $r=n$.  Thus, $H=\S_n$.  

\pfitem{\ref{FG2}}  If not, then we may assume by symmetry that $n(g\phi)=n$ for all $g\in G$.  Since $\im(\phi)\sub\S_n$ by~\ref{FG1}, it follows that $G$ embeds in $\S_{n-1}$, contradicting $n=\mu(G)$.

\pfitem{\ref{FG3}}  Let $\ze\in\T_{n+1}$ be the constant map with image $\{n+1\}$.  Then $\psi:S\to\T_{n+1}$, defined by $0\psi=\ze$ and $g\psi=g\phi\op\id_{\{n+1\}}$ for $g\in G$, is an embedding, so that $\mu(S)\leq n+1$.

To show that $\mu(S)\geq n+1$, suppose to the contrary that there is an embedding $\th:S\to\T_n$.  Let $\ze=0\th$, and let $i\in\im(\ze)$ be arbitrary, so $i=i\ze$.  By \ref{FG2} there exists $g\in G$ such that $i(g\th)\not=i$, so it follows that $i = i\ze = i(0\th) = i(0\th)(g\th) = i(g\th) \not=i$, a contradiction.
\epf

\begin{eg}
Let $G$ be a cyclic group of prime order $p$, and let $S=G\cup\{0\}$ be the semigroup obtained by adjoining a zero element to $G$.  It is well known \cite{Johnson1971} that $\mu(G)=p$, so it follows from Lemma \ref{lem:FG} that $\mu(S)=p+1$.
On the other hand, the variant $S^0$ (with sandwich element~$0$) is a null semigroup of size $p+1$.  Thus, for $p\geq5$ we have $\mu(S^0)<p+1=\mu(S)$, as follows quickly from Theorem~\ref{thm:null2}; cf.~Table \ref{tab:xial}.  

In fact, the rapid growth of the $\xi$ function means that the ratio $\mu(S)/\mu(S^0)$ can be made arbitrarily large by choosing $p$ large enough.  For example, again consulting Table \ref{tab:xial}, if $p$ has $20$ digits, then so too does $\mu(S)=p+1$, while $\mu(S^0)=30$.  By contrast, we have $\mu(S^a)/\mu(S)<2$ for any semigroup $S$ and any $a\in S$, as noted earlier.
\end{eg}

We can also use null and nilpotent semigroups to provide an example of a semigroup whose degree is strictly greater than the degrees of \emph{all} its variants.

\begin{eg}
Consider the $3$-nilpotent semigroup $S=\{x,y,z,0\}$ for which the only non-zero product is $x^2=y$:
\[
\begin{array}{c|cccc}
\cdot & 0&x&y&z\\
\hline
0&0&0&0&0\\
x&0&y&0&0\\
y&0&0&0&0\\
z&0&0&0&0\\
\end{array}
\]
We claim that:
\ben
\item \label{5} $\mu(S)=5$, but
\item \label{4} $\mu(S^a)=4$ for all $a\in S$.
\een
We first note that \ref{4} is clear, given the fact that any variant of a $3$-nilpotent semigroup is null, and that $\xi(3)=2$ and $\xi(4)=4$; cf.~Theorem \ref{thm:null2} and Table \ref{tab:xial}.

For \ref{5}, we have $\mu(S)\leq5$ because we have an embedding $S\to\T_5$ given by
\[
x \mt \trans{1&2&3&4&5\\1&1&1&2&4} \COMMA
y \mt \trans{1&2&3&4&5\\1&1&1&1&2} \COMMA
z \mt \trans{1&2&3&4&5\\1&1&2&1&1} \COMMA
0 \mt \trans{1&2&3&4&5\\1&1&1&1&1} .
\]
(If we denote these transformations by $x'$, $y'$, $z'$ and $0'$, respectively, then the only product not equal to $0'$ is $x'\circ x'=y'$.)  

To show that $\mu(S)\geq5$, we must show that $\T_4$ contains no subsemigroup isomorphic to $S$.  To do so, and aiming for a contradiction, suppose $\T_4$ does contain such a subsemigroup $T\cong S$.  Since $\{y,z,0\}$ is a null subsemigroup of $S$, it follows that $T$ contains a null subsemigroup of size~$3$.  For $2\leq r\leq 4$, the biggest null subsemigroup of $\T_4$ with a zero of rank $r$ has size $\nu_r(4)\leq\nu_2(4)=\xi(3)=2$; cf.~Theorem \ref{thm:null}.  It follows that the zero $\ze$ of $T$ has rank~$1$.  By symmetry, we may assume that $\ze=\trans{\bfour\\1}$.  As in Subsection \ref{subsect:maxnil}, it follows that $T$ is contained in the subsemigroup
\[
N = N(V,W) = \bigset{f\in\T_4}{\bfour f\sub W,\ Wf\sub V,\ Vf=\{1\}} \qquad\text{for some $\{1\}\subset V\subset W\subset \bfour$.}
\]
Again, by symmetry we can assume that $V=\btwo$ and $W=\bthree$.  Denoting the elements of $N$ by
\[
\ze = \trans{1&2&3&4\\1&1&1&1} \COMMa
a = \trans{1&2&3&4\\1&1&1&2} \COMMa
b = \trans{1&2&3&4\\1&1&1&3} \COMMa
c = \trans{1&2&3&4\\1&1&2&1} \COMMa
d = \trans{1&2&3&4\\1&1&2&2} \COMMa
e = \trans{1&2&3&4\\1&1&2&3} ,
\]
the multiplication table for $N$ is as follows:
\[
\begin{array}{c|cccccc}
\cdot & \ze&a&b&c&d&e\\
\hline
\ze&\ze&\ze&\ze&\ze&\ze&\ze\\
a&\ze&\ze&\ze&\ze&\ze&\ze\\
b&\ze&\ze&\ze&a&a&a\\
c&\ze&\ze&\ze&\ze&\ze&\ze\\
d&\ze&\ze&\ze&\ze&\ze&\ze\\
e&\ze&\ze&\ze&a&a&a\\
\end{array}
\]
Since $S$ contains an element with non-zero square, it follows that $e\in T$.  Since all other products in $T$ must equal $\ze$, it follows from the last row of the table that $T=\{\ze,a,b,e\}$.  But this $T$ has precisely two non-zero products, $b\circ e=e\circ e=a\not=\ze$.  This is the desired contradiction
\end{eg}

\section{\boldmath The degree of $\T_n^a$}\label{sect:Tna}

We now turn to variants of finite full transformation semigroups.  Recall that for a fixed transformation $a\in\T_n$, the variant $\T_n^a$ is the semigroup with underlying set $\T_n$ and operation $\star$ defined by $g\star h=gah$ for all $g,h\in\T_n$.  Recall from \eqref{eq:2n-r} that when $\rank(a)=r$, we have 
\[
n\leq\mu(\T_n^a)\leq 2n-r.
\]
The main guiding theme of the current section is Question 4.2 of \cite{JE2020}, which asks whether the upper bound of $2n-r$ is in fact the exact value of $\mu(\T_n^a)$.  

We begin in Subsection \ref{subsect:general} with some general results showing that any embedding of $\T_n^a$ in~$\T_{2n-r-1}$ (which would `break' the upper bound of $2n-r$) is rather restricted.  We then apply these results in Subsection \ref{subsect:r>1} to show that indeed $\mu(\T_n^a)=2n-r$ when $r$ is suitably large ($r\geq n-6$); see Theorem \ref{thm:n-6}.  

Subsection~\ref{subsect:r=1} then considers the other extreme case of $r=1$, where $a$ is a constant map.  Since $\mu(\T_n^a)$ does not depend on the actual choice of the constant map $a$, we denote it by $\mu(n)$, and the above bounds become $n\leq\mu(n)\leq2n-1$.  By realising $\T_n^a$ as an $n\times n^{n-1}$ right null semigroup (for ${\rank(a)=1}$), we apply Theorem \ref{thm:RN} to give a formula for $\mu(n)$ in Proposition~\ref{prop:mun}.  We then show in Proposition~\ref{prop:r=1} that it is indeed possible to break the upper bound of $2n-1$, and in Theorem~\ref{thm:mun} that the sequence $\mu(n)$ is strictly increasing in $n$.  Theorems \ref{thm:r=1xn} and \ref{thm:r=1n+k}, concern the asymptotic behaviour of $\mu(n)$, showing in particular that the ratio $\mu(n)/n$ tends to~$1$ as $n\to\infty$, while the difference $\mu(n)-n$ grows without bound.  We state some open problems in Subsection~\ref{subsect:OP}.

\subsection{General results}\label{subsect:general}

For the duration of this subsection, we fix an integer $n\geq2$, and a transformation $a\in\T_n$ with $\rank(a)=r$.  Since we wish to study the variant $\T_n^a$, and since $\T_n^a\cong\T_n^{gah}$ for any permutations $g,h\in\S_n$ (see \cite[Proposition~13.1.3]{GMbook}), we may assume without loss of generality that $a$ is an idempotent with $\im(a)=\br$.

As noted above, we have $\mu(\T_n^a)\leq 2n-r$, and here we are interested in the question of whether we could in fact have strict inequality: $\mu(\T_n^a)<2n-r$.  If we do, then $\T_n^a$ can be embedded in $\T_{2n-r-1}$.  For the duration of this subsection, we write $m=2n-r-1$ for simplicity, and we hypothesise the existence of an embedding
\[
\phi:\T_n^a\to\T_m.
\]
Since $\phi$ obviously cannot exist for $r=n$ (since then $m=n-1$), we assume that $r<n$ in all that follows.
We will write $f'=f\phi$ for all $f\in\T_n$, noting that $(f\star g)'=f'g'$ for all $f,g\in\T_n$.  If $S\sub\T_n$, we write $S'=S\phi=\set{f'}{f\in S}$.  

We will occasionally need to refer to Green's relations on the semigroups $\T_n^a$ and $\T_m$; those on the former will be denoted $\L^a$, $\R^a$, and so on, to distinguish them from those on the latter, which will be denoted in the usual way.

Consider the set
\[
S = \set{f\in\T_n}{f|_{\br}=\id_{\br}} = \set{f\in\T_n}{af=a},
\]
which is easily seen to be a subsemigroup of $\T_n^a$ of size $n^{n-r}$.  Consequently, $S'=S\phi$ is a subsemigroup of $\T_m$ of the same size.  Also define
\[
F = \bigcup_{f\in \T_n}\im(f') \AND G = \bigcup_{f\in S}\im(f'),
\]
noting that $G\sub F\sub\bm$.  The following is true by definition:

\begin{lemma}\label{lem:TF}
The image $\im(\phi)$ is contained in the subsemigroup
\[
\epfreseq
\T_m(F) = \set{f\in\T_m}{\im(f)\sub F}.
\]
\end{lemma}

\begin{lemma}\label{lem:G}
For all $f\in S$ and all $x\in F$ we have $xf'=xa'$.  Consequently, $n^{n-r} \leq |G|^{m-|F|}$.
\end{lemma}

\pf
Fix some $f\in S$ and $x\in F$.  By definition, we have $x=yg'$ for some $y\in\bm$ and $g\in\T_n$.  Since $f\in S$ we have $a=af$, so $g\star f=gaf=ga=gaa=g\star a$, which gives $g'f'=g'a'$ in $\T_m$.  But then $xf' = yg'f' = yg'a' = xa'$.

Now that we have proved the first assertion, it follows that the elements of $S'$ are only distinguished by their restriction to $\bm\sm F$.  Every such restriction is a map $\bm\sm F\to G$, of which there are $|G|^{m-|F|}$.  Since $|S'|=n^{n-r}$, the second assertion follows.
\epf

\begin{lemma}\label{lem:F}
We have $|F|\leq n-2$.  
\end{lemma}

\pf
Suppose to the contrary that $|F|=n-1+k$ where $k\geq0$.  Then by Lemma \ref{lem:G} and Corollary~\ref{cor:fact1}, and since $G\sub F$, we have
\[
n^{n-r} \leq |G|^{m- |F|} \leq |F|^{m-|F|} = (n-1+k)^{(2n-r-1)-(n-1+k)} = (n-1+k)^{n-r-k} \leq (n-1)^{n-r},
\]
a contradiction (as $n-r\geq1$).
\epf

In what follows, a crucial role will be played by the constant maps from $\T_n$, and their images in $\T_m$ under $\phi$.
For $i\in\bn$, write $e_i=\binom{\bn}i\in\T_n$ for the constant map with image $\{i\}$.  We will frequently make use of the fact that $e_if=e_{if}$ for all $i\in\bn$ and $f\in\T_n$.  

Since the $e_i$ are all $\R^a$-related in $\T_n^a$, the $e_i'=e_i\phi$ are all $\R$-related in $\T_m$, so they have a common kernel, the classes of which we will denote by $B_1,\ldots,B_q$.  Thus, we may write
\[
e_i' = \trans{B_1&\cdots&B_q\\x_{i1}&\cdots&x_{iq}} \qquad\text{for $i\in\bn$.}
\]
Since each $e_i'$ is an idempotent, we have $x_{ij}\in B_j$ for all $i\in\bn$ and $j\in\bq$.

\begin{lemma}\label{lem:q}
We have $q\geq2$.
\end{lemma}

\pf
If we had $q=1$ then $e_1',\ldots,e_n'$ would be distinct transformations of rank $1$, and so $|F|\geq|\im(e_1')\cup\cdots\cup\im(e_n')|=n$, contradicting Lemma \ref{lem:F}.
\epf

The next statement uses the $\op$ operation defined in Subsection \ref{subsect:TX}.

\begin{lemma}\label{lem:B_i}
For any $f\in\T_n$ and any $i\in\bq$, we have $B_if'\sub B_i$.  Consequently, $\im(\phi)$ is contained in $\T_{B_1}\op\cdots\op\T_{B_q}$.
\end{lemma}

\pf
In $\T_n^a$ we have $f\star e_1=e_1$, so it follows that $f'e_1'=e_1'$.  The claim follows.
\epf

For each $i\in\bq$ let $F_i = F\cap B_i$.  The next result follows immediately from Lemmas \ref{lem:TF} and~\ref{lem:B_i}.

\begin{cor}\label{cor:B_i}
The image $\im(\phi)$ is contained in $\T_{B_1}(F_1)\op\cdots\op\T_{B_q}(F_q)$.  \epfres
\end{cor}

Next we wish to find a lower bound on the rank of $a'$.  For this, we require the following result concerning Green's $\leqJ$-preorder, defined on a semigroup $S$ by $x\leqJ y \iff x\in S^1yS^1$.  In a full transformation semigroup $\T_X$, we have $f\leqJ g \iff\rank(f)\leq\rank(g)$.  The proof uses the natural partial order on the idempotents $E(S)$ defined by $e\leq f \iff e=ef=fe$.

\begin{lemma}\label{lem:EH}
Suppose we have an embedding $\psi:S\to T$, where $S$ and $T$ are finite regular semigroups, and where all elements of $S$ are comparable in the $\leqJ$-preorder.  Then for all ${x,y\in S}$, we have $x\leqJ y$ in $S$ if and only if $x\psi \leqJ y\psi$ in $T$.  
\end{lemma}

\pf
Clearly $x\leqJ y\implies x\psi\leqJ y\psi$.  For the converse, suppose $x\psi\leqJ y\psi$.  By assumption we have either $x\leqJ y$ or $x\geqJ y$.  In the former case we are done, so suppose instead that $x\geqJ y$.  This implies that $x\psi\geqJ y\psi$, so in fact $x\psi\J y\psi$.  Since $x\geqJ y$ in the regular semigroup $S$, it follows from \cite[Theorem 1]{Hall1970} that there exist idempotents $e,f\in E(S)$ such that $x\J e\geq f\J y$.  This implies that the idempotents $e\psi,f\psi\in E(T)$ satisfy $e\psi\geq f\psi$, and also $e\psi\J x\psi \J y\psi\J f\psi$.  As $T$ is finite, the natural order is trivial on $\J$-classes (see for example \cite[Result 6]{Hall1973}), so that $e\psi=f\psi$, and thus $e=f$ as $\psi$ is injective.  We thus have $x\J e=f\J y$.  In particular, $x\leqJ y$, as required.
\epf

\begin{lemma}\label{lem:q+r-1}
We have $\rank(a')\geq q+r-1$.
\end{lemma}

\pf
Consider the set $R=\Reg(\T_n^a)$ of all regular elements of $\T_n^a$.  As shown in \cite{DE2015},~$R$ is a subsemigroup of $\T_n^a$, and is a chain of ${\J^a}={\D^a}$-classes: $D_1^a<\cdots<D_r^a$.  Fix some $f_i\in D_i^a$ for each ${i\in\br}$, and assume that $f_1=e_1$ and $f_r=a$.  Then using Lemma \ref{lem:EH} and the above-mentioned characterisation of the $\leqJ$-ordering in $\T_m$,
\begin{align*}
f_1 <_{\J}\cdots<_{\J}f_r &\implies f_1'<_{\J}\cdots<_{\J}f_r' \\
&\implies q=\rank(f_1')<\rank(f_2')<\cdots<\rank(f_r')=\rank(a').
\end{align*}
The result is then immediate.
\epf

In the case $r=1$, $a$ is a constant map, so that $\rank(a')$ attains the lower bound of $q+r-1=q$ just established.  We show below that the lower bound is never attained when $r\geq2$; see Lemma~\ref{lem:q+r}.

First, however, we take a closer look at the transformation $a'=a\phi$.
For each $i\in\bq$, denote the restriction $a'|_{B_i}$ by $a_i'$, and write $r_i=\rank(a_i')$.  By Lemma~\ref{lem:B_i} we have $a' = a_1'\op\cdots\op a_q'$, with each $a_i'\in\T_{B_i}(F_i)$.  

\begin{lemma}\label{lem:r1r}
We have $r_i\geq r$ for some $i\in\bq$.
\end{lemma}

\pf
Recall that $\phi$ is an embedding of $\T_n^a$ into $\T_{B_1}\op\cdots\op\T_{B_q}$.  For $i\in\bq$, let $\phi_i:\T_n^a\to\T_{B_i}$ be the result of composing $\phi$ with the projection onto the $i$th coordinate.  As in \cite[Section~4.3]{Sandwiches2}, the set $T=a\T_na=a\star\T_n\star a$ is a subsemigroup of $\T_n^a$ isomorphic to $\T_r$.  Consider the restrictions $\psi=\phi|_T$ and $\psi_i=\phi_i|_T$ for each $i\in\bq$, noting that $\psi$ is injective.  If all the $\psi_i$ were non-injective, then each $\ker(\psi_i)$ would contain the minimum non-trivial congruence on $T$ (the congruences on $T\cong\T_r$ form a chain \cite{Malcev1952}), but then so too would $\ker(\psi)=\ker(\psi_1)\cap\cdots\cap\ker(\psi_q)$, meaning that $\psi$ is non-injective, a contradiction.  It follows that some $\psi_i:T\to\T_{B_i}$ is injective.  If we write~$p$ for the minimum rank of an element in $\im(\psi_i)$, then as in the proof of Lemma \ref{lem:q+r-1}, we have $r_i=\rank(a_i')=\rank(a_i\psi)\geq p+r-1\geq r$.
\epf

Thus, without loss of generality we may assume that $r_1\geq r$.  It follows as well that $|F_1|\geq r$ (as $a_1'\in\T_{B_1}(F_1)$ and $\rank(a_1')\geq r$).

\begin{lemma}\label{lem:r_i}
If $r\geq2$, then $r_i\geq2$ for some $2\leq i\leq q$.
\end{lemma}

\pf
Seeking a contradiction, suppose instead that $r_2=\cdots=r_q=1$.  Then $B_2,\ldots,B_q$ are all $\ker(a')$-classes.  Let the $\ker(a_1')$-classes be $A_1,\ldots,A_{r_1}$.
Thus, we may write
\[
a' = \trans{A_1&\cdots&A_{r_1}&B_2&\cdots&B_q\\y_1&\cdots&y_{r_1}&z_2&\cdots&z_q},
\]
where each $y_i\in A_i$ ($1\leq i\leq r_1$) and $z_j\in B_j$ ($2\leq j\leq q$), as $a'$ is an idempotent.  Linking back to the $e_i'$ ($1\leq i\leq n$), we claim that
\begin{equation}\label{eq:zj}
x_{ij}=z_j \qquad\text{for all $1\leq i\leq n$ and $2\leq j\leq q$.}
\end{equation}
To prove this, fix some $1\leq i\leq n$ and $2\leq j\leq q$.

Suppose first that $1\leq i\leq r$.  
Since $a$ is an idempotent with image $\br$ we have $i=ia$, and so $e_i=e_ia=e_iaa=e_i\star a$.  Thus, $e_i'=e_i'a'$, and so $x_{ij}\in\im(e_i')=\im(e_i'a')\sub\im(a')$.
But $x_{ij}\in B_j$, and since $\im(a')\cap B_j=\{z_j\}$ it follows that $x_{ij}=z_j$.

Now suppose $r<i\leq n$.  Put $k=ia$, and let $l\in\br\sm\{k\}$ be arbitrary (recall that $r\geq2$).  Choose any $f\in\T_n$ such that $kf=i$ and $lf=l$.  Since $k,l\in\im(a)$, and since $a$ is an idempotent, we have $k=ka$ and $l=la$.  It follows that $e_k\star f = e_i$ and $e_l\star f=e_l$, so that $e_k'f'=e_i'$ and $e_l'f'=e_l'$.  It follows from these, respectively, that $x_{kj}f'=x_{ij}$ and $x_{lj}f'=x_{lj}$.  But $k,l\leq r$, so by the previous paragraph we have $x_{kj}=z_j=x_{lj}$, and it follows that $x_{ij}=x_{kj}f'=x_{lj}f'=x_{lj}=z_j$.  This completes the proof of \eqref{eq:zj}.

It now follows that
\[
e_i' = \trans{B_1&B_2&\cdots&B_q\\x_{i1}&z_2&\cdots&z_q} \qquad\text{for all $i\in\bn$.}
\]
Since the $e_i'$ are pairwise distinct, it follows that the set $\{x_{11},\ldots,x_{n1}\}$ has size $n$.  But this set is contained in $F$, and this contradicts Lemma \ref{lem:F}.
\epf

\begin{lemma}\label{lem:q+r}
If $r\geq2$, then $\rank(a')\geq q+r$.
\end{lemma}

\pf
We have $\rank(a')\geq q+r-1$ by Lemma \ref{lem:q+r-1}.  If in fact $\rank(a')=q+r-1$, then from 
\[
q+r-1 = \rank(a') = r_1+r_2+\cdots+r_q \geq r+r_2+\cdots+r_q ,
\]
it follows that $r_2+\cdots+r_q\leq q-1$, so that $r_2=\cdots=r_q=1$, contradicting Lemma \ref{lem:r_i}.
\epf

The following is a simple consequence of Lemmas \ref{lem:q} and \ref{lem:q+r}:

\begin{cor}\label{cor:r+2}
If $r\geq2$, then $\rank(a')\geq r+2$.  \epfres
\end{cor}

Now that we have gathered all the general results we need, we proceed to consider separate cases for `large' and `small' $r$ in the next two subsections.

\subsection[The case $r\geq n-6$]{\boldmath The case $r\geq n-6$}\label{subsect:r>1}

We now consider the case in which the rank $r$ of the sandwich element $a$ is suitably large, meaning specifically that $r\geq n-6$.  We will shortly prove that for such large $r$, the degree $\mu(\T_n^a)$ reaches its upper bound of $2n-r$.  For this we need the following simple lemma.

\begin{lemma}\label{lem:lam}
For any composition $\si=(s_1,\ldots,s_r)\vDash n$ we have $\prod\si\geq n-r+1$.
\end{lemma}

\pf
First note that for any $s\geq t\geq2$ we have $st\geq2s>2s-1\geq s+t-1=(s+t-1)\cdot1$.  It quickly follows that $\prod\si\geq\prod\si'$ for the composition $\si'=(n-r+1,1,\ldots,1)$.
\epf

\begin{thm}\label{thm:n-6}
If $a\in\T_n$ with $r=\rank(a)\geq n-6$, then $\mu(\T_n^a)=2n-r$.
\end{thm}

\pf
Again we write $m=2n-r-1$, and aiming for a contradiction, we assume there is an embedding $\phi:\T_n^a\to\T_m$.  As usual, we may assume that $a$ is an idempotent, and we keep the notation of Subsection \ref{subsect:general}, writing $f'=f\phi$, and so on.  Keeping Lemmas \ref{lem:TF} and \ref{lem:F} in mind, we may assume by symmetry that $\phi$ maps $\T_n^a$ into $\T_m(n-2)$.  

Let $t=\rank(a')$, so that $r+2\leq t\leq n-2$, by Corollary \ref{cor:r+2}.  We immediately obtain a contradiction if $r+2>n-2$: i.e., if $r\geq n-3$.  So for the rest of the proof we assume that $n-6\leq r\leq n-4$.  Incidentally, this implies $n\geq5$.

Let $D$ be the $\D^a$-class of $a$ in $\T_n^a$.  Since all elements of $D'=D\phi$ are $\D$-related (in $\T_m$) to~$a'$, and since $\rank(a')=t$, it follows that 
\[
D' \sub D_t(\T_m(n-2)) = \set{f\in\T_m(n-2)}{\rank(f)=t}.
\]
(Note that $D_t(\T_m(n-2))$ is not itself a $\D$-class of $\T_m(n-2)$.)
By \cite[Theorem~5.7(v)]{DE2015}, $D$ is an $r^{n-r}\times\Lam$ rectangular group over $\S_r$, where $\Lam$ is the product of the sizes of the $\ker(a)$-classes.  By Lemma~\ref{lem:lam} we have $\Lam\geq n-r+1$, so it follows that $D$ contains an $r^{n-r}\times(n-r+1)$ rectangular band.  So too therefore does $D'$, and hence also $D_t(\T_m(n-2))$.  It follows, using the~$\lam$ and $\rho$ parameters from Subsection \ref{subsect:LZRZ}, that
\begin{align*}
\lam_t(m,n-2) \geq r^{n-r} &\AND \rho_t(m,n-2)\geq n-r+1 . 
\intertext{Using Lemma \ref{lem:Tnk}, this is equivalent to}
\lam_t(m) \geq r^{n-r} &\AND \rho_t(n-2)\geq n-r+1 . 
\intertext{For the rest of the proof we write $r=n-k$, noting that $k\in\{4,5,6\}$.  The previous inequalities become}
\lam_t(n+k-1) \geq (n-k)^k &\AND \rho_t(n-2)\geq k+1 . 
\end{align*}
We will obtain the desired contradiction by showing that 
\begin{equation}\label{eq:rho}
\rho_t(n-2)< k+1 \qquad\text{for all $n-k+2=r+2\leq t\leq n-2$.}
\end{equation}
To do so, we consider the allowable values of $k$ separately.  In the following we make use of Remark \ref{rem:GM1'}, which tells us that $\rho_{n-2-l}(n-2)=2^l$ if $0\leq l\leq\frac{n-2}2$.  We only wish to apply this for $l\in\{0,1,2\}$, so since $n\geq5$, this could only be invalid when $n=5$ and $l=2$; however, in this case we have $\rho_{n-2-l}(n-2)=\rho_1(3)=3<4=2^l$, which is sufficient for our purposes.  
\bit
\item When $k=4$, we only have $t=n-2$ to consider, and \eqref{eq:rho} holds since $\rho_{n-2}(n-2)=1 < 5$.
\item When $k=5$, we only have $t=n-3$ and $n-2$ to consider, and \eqref{eq:rho} holds since
\[
\rho_{n-3}(n-2) = 2<6 \AND \rho_{n-2}(n-2)=1 < 6.
\]
\item When $k=6$, \eqref{eq:rho} holds since
\[
\rho_{n-4}(n-2) \leq 4<7 \COMMA \rho_{n-3}(n-2) = 2<7  \COMMA \rho_{n-2}(n-2)=1 < 7.
\]
(The $\leq$ sign in the first of these is because of the $(n,l)=(5,2)$ case mentioned above.)  \qedhere
\eit
\epf

\begin{rem}
The argument in the above proof breaks down when $r=n-7$.  Here $k=7$, and~\eqref{eq:rho} does not hold for the minimum value of $t=n-5$, as we have $\rho_{n-5}(n-2) = 8 = k+1$.  We do not currently know if Theorem \ref{thm:n-6} holds for $r=n-7$, or more generally.
\end{rem}

\subsection[The case $r=1$]{\boldmath The case $r=1$}\label{subsect:r=1}

The previous subsection considered the case in which $r\geq n-6$ was suitably large.  We now consider the other extreme case in which $r=1$, meaning that the sandwich element $a$ is a constant map.  Since the value of $\mu(\T_n^a)$ does not depend on the particular choice of the constant map $a$, we denote it by $\mu(n)=\mu(\T_n^a)$.  Here \eqref{eq:2n-r} gives the bounds $n\leq\mu(n)\leq 2n-1$.

The variant $\T_n^a$ has a rather simple structure when $\rank(a)=1$.  To describe it, we may assume by symmetry that $a=e_1=\binom{\bn}1$.  For each $i\in\bn$, define
\[
S_i = \set{f\in\T_n}{1f=i}.
\]
Then $\T_n=S_1\cup\cdots\cup S_n$, with $|S_i|=n^{n-1}$, and we note that only $S_1$ is a subsemigroup of $\T_n$.  However, we have
\begin{equation}\label{eq:SiSj}
S_i\star S_j = \{e_j\} \qquad\text{for all $i,j\in\bn$,}
\end{equation}
and this determines the entire multiplication table of $\T_n^a$.  In particular, each $S_i$ is a null subsemigroup of $\T_n^a$ with zero $e_i$.  Moreover, \eqref{eq:SiSj} exhibits $\T_n^a$ as an $n\times n^{n-1}$ (uniform) right null semigroup, in the language of Subsection \ref{subsect:RN}.  It therefore follows that $\mu(n)=\mu(\T_n^a)=\rn(n,n^{n-1})$, so we have the following:

\begin{prop}\label{prop:mun}
For $n\geq2$ we have
\begin{align}
\label{eq:mun1} \mu(n) &= \min \set{2A + 3B + C + D }{ 0\leq A\leq 2,\ B,D\geq0,\ C\geq2,\ 2^A3^BC\geq n,\ C^D\geq n^{n-1} }\\
\label{eq:mun2} &= {\min} \bigset{2A + 3\lceil\log_3(n/2^AC)\rceil + C + \lceil\log_Cn^{n-1}\rceil }{ 0\leq A\leq 2\leq  C \leq2n-1},
\end{align}
where the parameters $A,B,C,D$ in the above sets are all integers.
\end{prop}

\pf
This follows immediately from Theorem \ref{thm:RN}\ref{RN2} and Remark \ref{rem:RN}.  The latter explains that the upper bound on $C$ in \eqref{eq:mun2} can be any known upper bound on $\rn(n,n^{n-1})=\mu(n)$, and we know from Theorem \ref{thm:EM} that $\mu(n)\leq2n-1$.  
\epf

Table \ref{tab:mun} gives some values of $\mu(n)$, computed using \eqref{eq:mun2}.

\begin{table}[H]
\begin{center}
\scalebox{0.97}{
\begin{tabular}{|c|cccccccccccccccccccc|}
\hline
$n$ & 1 & 2 & 3 & 4 & 5 & 6 & 7 & 8 & 9 & 10 & 11 & 12 & 13 & 14 & 15 & 16 & 17 & 18 & 19 & 20\\
\hline
$\mu(n)$ & 1 & 3 & 5 & 7 & 9 & 11 & 13 & 15 & 17 & 19 & 21 & 23 & 25 & 27 & 29 & 30 & 32 & 34 & 36 & 37\\
\hline
\end{tabular}
}
\\[5mm]
\begin{tabular}{|c|rrrrrrrrrr|}
\hline
$n$ & \phantom{1}100 & \phantom{1}200 & \phantom{1}300 & \phantom{1}400 & \phantom{1}500 & \phantom{1}600 & \phantom{1}700 & \phantom{1}800 & \phantom{1}900 & 1000 \\
\hline
$\mu(n)$ & \phantom{1}167 & \phantom{1}321 & \phantom{1}473 & \phantom{1}623 & \phantom{1}772 & \phantom{1}919 & 1066 & \phantom{1}1213 & \phantom{1}1359 & \phantom{1}1504 \\
\hline
\end{tabular}
\\[5mm]
\begin{tabular}{|c|rrrrrrrrrr|}
\hline
$n$ & 1100 & 1200 & 1300 & 1400 & 1500 & 1600 & 1700 & \phantom{1}1800 & \phantom{1}1900 & \phantom{1}2000 \\
\hline
$\mu(n)$ &1649 & 1794 & 1938 & 2082 & 2226 & 2369 & 2513 & \phantom{1}2656 & \phantom{1}2798 & \phantom{1}2941 \\
\hline
\end{tabular}
\\[5mm]
\begin{tabular}{|c|cccccccccc|}
\hline
$n$ & 1000 & 2000 & 3000 & 4000 & 5000 & 6000 & 7000 & \phantom{1}8000 & \phantom{1}9000 & 10000 \\
\hline
$\mu(n)$ &1504 & 2941 & 4359 & 5766 & 7165 & 8557 & 9945 & 11329 & 12709 & 14086 \\
\hline
\end{tabular}
\caption{Values of $\mu(n)$, which is the degree of a variant $\T_n^a$ for any $a\in\T_n$ of rank $1$.}
\label{tab:mun}
\end{center}
\end{table}

In what follows, we will sometimes establish an upper bound $\mu(n)\leq m$ by showing that $m=2A+3B+C+D$ for some tuple $(A,B,C,D)$ satisfying the restrictions listed in \eqref{eq:mun1}.  We will say that such a tuple \emph{witnesses} the inequality $\mu(n)\leq m$.

\begin{rem}\label{rem:improv}
We can improve \eqref{eq:mun2} to:
\[
\mu(n) = {\min} \bigset{2A + 3\lceil\log_3(n/2^AC)\rceil + C + \lceil\log_Cn^{n-1}\rceil }{ 0\leq A\leq 2\leq  C \leq n},
\]
where the bound $C\leq2n-1$ is replaced by $C\leq n$.  Indeed, suppose there is a witness $(A,B,C,D)$ to $\mu(n)=m$ with $C>n$.  By minimality, and since $2^03^0C>n$, we must have ${A=B=0}$.  
Write $t=C-n>0$ and let $k=D+t\geq0$.  Then ${2n-1\geq \mu(n)=C+D=n+t+D=n+k}$, so we obtain $k\leq n-1<n$.  Corollary \ref{cor:fact1} then gives 
\[
C^D = (n+t)^{k-t} \leq n^k,
\]
and since also $n+k=C+D=\mu(n)$, it follows that $(0,0,n,k)$ also witnesses $\mu(n)=m$.
\end{rem}

We are now in a position to show that the upper bound $\mu(\T_n^a)\leq2n-r$ coming from Theorem~\ref{thm:EM} is not necessarily the exact value.

\begin{prop}\label{prop:r=1}
\ben
\item \label{r=11} If $n\leq15$, then $\mu(n)=2n-1$.
\item \label{r=12} If $n\geq16$, then $\mu(n)\leq2n-2$.
\een
\end{prop}

\pf
\firstpfitem{\ref{r=11}}  This follows by inspecting Table \ref{tab:mun}.

\pfitem{\ref{r=12}}  By Table \ref{tab:mun}, it suffices to prove this for $n\geq19$.  To do this, let $k=\lceil\frac n2\rceil$.  We will show that the tuple $(A,B,C,D) = (1,0,k,2n-4-k)$ witnesses the inequality $\mu(n)\leq2n-2$.  Since $2A+3B+C+D=2n-2$, it remains to show that $A,B,C,D$ satisfy the restrictions listed in~\eqref{eq:mun1}.  These are all completely routine, apart from $C^D\geq n^{n-1}$, for which we have
\begin{align*}
C^D = k^{2n-4-k} &\geq \left(\frac n2\right)^{2n-4-\frac{n+1}2} &&\text{as $\tfrac n2\leq k\leq\tfrac{n+1}2$}\\
&= n^{n-1} \cdot \frac{n^{(n-7)/2}}{2^{(3n-9)/2}} \\
&\geq n^{n-1} \cdot \frac{2^{2n-14}}{2^{(3n-9)/2}} &&\text{as $n\geq16=2^4$}\\
&= n^{n-1} \cdot 2^{(n-19)/2} \\
&\geq n^{n-1} &&\text{as $n\geq19$.}  \qedhere
\end{align*}
\epf

The upper bound of $2n-2$ in Proposition \ref{prop:r=1}\ref{r=12} is still not sharp in general, as clearly indicated in Table \ref{tab:mun}.  In fact, given any integer~$k$, we have $\mu(n)\leq 2n-k$ for suitably large~$n$; indeed, in Theorem \ref{thm:r=1xn} below, we prove an even stronger statement.
First, however, we demonstrate another interesting property of the numbers $\mu(n)$, namely that they are strictly increasing:

\begin{thm}\label{thm:mun}
We have $\mu(1)<\mu(2)<\mu(3)<\cdots$.
\end{thm}

\pf
Let $n\geq2$, and write $m=\mu(n)$.  We must show that $\mu(n-1)\leq m-1$.  Consulting Table~\ref{tab:mun}, it suffices to assume that $n\geq16$, which we do for the rest of the proof.  By Proposition~\ref{prop:r=1} we have $m\leq2n-2$.

Suppose $\mu(n)=m$ is witnessed by the tuple $(A,B,C,D)$; as in Remark \ref{rem:improv}, we can assume that $C\leq n$.  We aim to show that the inequality $\mu(n-1)\leq m-1$ is witnessed by $(A,B,C,D-1)$.  Again, verification of the required conditions is mostly routine; the only exceptions this time are $D-1\geq0$ and $C^{D-1}\geq (n-1)^{n-2}$.


Keeping $C\leq n$ in mind, we have $n^{n-1}\leq C^D\leq n^D$.  It follows that $D\geq n-1$, so certainly $D-1\geq0$.  We also have 
\[
n\cdot n^{n-2} \leq C^D = C\cdot C^{D-1} \leq n\cdot C^{D-1},
\]
which implies $C^{D-1}\geq n^{n-2} \geq (n-1)^{n-2}$.  As noted above, this completes the proof.
\epf

The next two theorems concern the asymptotic behaviour of the sequence $\mu(n)$.  For their proofs, it will be convenient to use the following special case of Lemma \ref{lem:mpq}.

\begin{lemma}\label{lem:stm}
For $m,n\in\N$ we have $\mu(n)\leq m$ if and only if there exist compositions $\si$ and $\tau$ such that
\begin{align}
\label{eq:stm}
\tau\pre\si\vDash m \COMMA \textstyle{\prod\tau\geq n} \COMMA \eta(\si,\tau)\geq n^{n-1}.\\
\epfreseq
\end{align}
\end{lemma}


\begin{thm}\label{thm:r=1xn}
For any real $1<x<2$, there exists $N\in\N$ such that $\mu(n)\leq xn$ for all $n>N$.  Thus, asympotically we have
\[
\mu(n) \sim n \qquad\text{as $n\to\infty$.}
\]
\end{thm}

\pf
By Lemma \ref{lem:stm}, it suffices to show that there exists $N\in\N$ such that for all $n>N$ there exist compositions
\[
\tau_n\pre\si_n\vDash  \lfloor xn\rfloor \qquad\text{such that}\qquad \textstyle{\prod\tau_n}\geq n \AND \eta(\si_n,\tau_n)\geq n^{n-1}.
\]

Let $0<y<x-1$ be arbitrary, and set $z=x-y$, noting that $z>1$.  Also let $k\in\N$ be such that $ky\geq1$.  Note that $y$, $z$ and $k$ depend only on $x$.

For any $n>k+2$, define
\[
\si_n = (\lfloor xn\rfloor-k,k) \AND \tau_n = (\lceil yn\rceil,k).
\]
Clearly $\prod\tau_n\geq (ky)n\geq n$.  Since
\[
\lfloor xn\rfloor -k > (xn-1) - k = n +(x-1)n  - 1 - k > (k+2) +yn -1-k = yn+1 > \lceil yn\rceil,
\]
it follows that $\si_n$ is a composition of $\lfloor xn\rfloor$, and that $\tau_n\pre\si_n$.  Next, recalling that $z=x-y>1$, we have
\[
\eta(\si_n,\tau_n) = \lceil yn\rceil ^{\lfloor xn\rfloor-k-\lceil yn\rceil} \geq (yn) ^{(xn-1)-k-(yn+1)} = (yn)^{zn-k-2}.
\]
Since $k,y>0$ and $z>1$ are constants (not depending on $n$), the function $(yn)^{zn-k-2}$ dominates~$n^{n-1}$, so there exists $N'\in\N$ such that $(yn)^{zn-k-2}>n^{n-1}$ for all $n>N'$.  We then take ${N=\max(k+2,N')}$.
\epf

Theorem \ref{thm:r=1xn} tells us that the ratio~$\mu(n)/n$ tends to $1$ as $n$ increases.  In contrast to this, the difference $\mu(n)-n$ grows without bound:

\begin{thm}\label{thm:r=1n+k}
For any integer $k\geq2$ we have $\mu(n)\geq n+k$ for all $n\geq\lceil k\ln k\rceil+1$.
\end{thm}

For the proof we require the following lemma, concerning the $L$ function defined in Subsection~\ref{subsect:numbers}:

\begin{lemma}\label{lem:L(m)}
If $\mu(n)\leq m$ where $n\geq2$, then $n^{n-1} < \lceil L(m)\rceil^{m-\lfloor L(m)\rfloor}$.
\end{lemma}

\pf
Fix compositions $\si$ and $\tau$ satisfying \eqref{eq:stm}.  Then using the definitions and Lemma \ref{lem:Xi}\ref{Xi1}, we have
\[
n^{n-1} \leq \eta(\si,\tau) \leq \xi(\si) \leq \Xi(m) = \xi(m).
\]
Now write $x=L(m)$, so that $\xi(m)=\max\left(\lfloor x\rfloor^{m-\lfloor x\rfloor},\lceil x\rceil^{m-\lceil x\rceil}\right)$ by Lemma \ref{lem:xi}\ref{xi3}.  Since $m$ is an integer, $x$ is not an integer, so it follows that $\lfloor x\rfloor < \lceil x\rceil$.  Since $m\geq\mu(n)>n\geq2$, we also have $m\geq x$, so that $m\geq\lceil x\rceil>\lfloor x\rfloor$.
All of this implies that $\lfloor x\rfloor^{m-\lfloor x\rfloor}$ and $\lceil x\rceil^{m-\lceil x\rceil}$ are both strictly less than $\lceil x\rceil^{m-\lfloor x\rfloor}$.
\epf

\pf[\bf Proof of Theorem \ref{thm:r=1n+k}]
Fix some $n\geq\lceil k\ln k\rceil+1$ where $k\geq2$, and let $x=L(n+k-1)$, so that $n+k-1=x(1+\ln x)$ by definition.  Note that $n\geq k+1$ since $k\geq2$.  

We first claim that $\lceil x\rceil\leq n$.
Indeed, this is clear if $x<e$ since then $n\geq k+1\geq3=\lceil e\rceil$.  Now suppose $x\geq e$, so that $1+\ln x\geq2$.  Then $2x\leq x(1+\ln x)=n+k-1\leq2n-1$, as $n\geq k$, and so $x<n$.  The claim follows as $n$ is an integer.

Next we claim that $\lfloor x\rfloor\geq k$.  For this we have
\[
x(1+\ln x) = n+k-1 \geq (\lceil k\ln k\rceil+1)+k-1 \geq k(\ln k+1).
\]
Since $x,k\geq e^{-1}$, it follows that $x\geq k$, and again the claim follows as $k$ is an integer.

Returning now to the main proof, suppose to the contrary that $\mu(n)\leq n+k-1$.  Then the above two claims and Lemma \ref{lem:L(m)} (with $m=n+k-1$) give
\[
n^{n-1} < \lceil x\rceil^{n+k-1-\lfloor x\rfloor} \leq n^{(n+k-1)-k} = n^{n-1},
\]
a contradiction.
\epf

\subsection{Open problems}\label{subsect:OP}

We conclude with a number of open problems.  The most obvious is the following:

\begin{prob}
Give a formula for $\mu(\T_n^a)$ for arbitrary $a\in\T_n$.  
\bit
\item Does $\mu(\T_n^a)$ depend only on $n$ and $\rank(a)$?  
\item Does $\mu(\T_n^a)=2n-r$ whenever $r=\rank(a)\geq2$?  
\item Classify the pairs $(n,r)$ for which $\mu(\T_n^a)=2n-r$ for all $a\in\T_n$ with $r=\rank(a)$.
\eit
\end{prob}

In the absence of a formula for $\mu(\T_n^a)$, it would also be interesting to answer the following:

\begin{prob}
Given $a,b\in\T_n$ with $r=\rank(a)$ and $s=\rank(b)$, which (if any) of the following implications hold?
\bit
\item $r<s\implies\mu(\T_n^a)<\mu(\T_n^b) $,
\item $r<s\implies\mu(\T_n^a)\leq\mu(\T_n^b) $,
\item $r\leq s\implies\mu(\T_n^a)\leq\mu(\T_n^b) $,
\item $r=s\implies\mu(\T_n^a)=\mu(\T_n^b) $.
\eit
\end{prob}

\begin{prob}
For which numbers $r\in\N$ (if any) do we have $\mu(\T_n^a)=2n-r$ for all $n\in\N$ and all $a\in D_r(\T_n)$?  For example, $r=1$ does not satisfy this property (Proposition \ref{prop:r=1}).
\end{prob}

\begin{prob}
For which numbers $k\in\N\cup\{0\}$ do we have $\mu(\T_n^a)=n+k$ for all $n\in\N$ and all $a\in D_{n-k}(\T_n)$?  For example, $k=0,1,\ldots,6$ all satisfy this property (Theorem \ref{thm:n-6}).
\end{prob}

\footnotesize
\def\bibspacing{-1.1pt}
\bibliography{biblio}
\bibliographystyle{abbrv}

\end{document}